\documentclass[a4paper,12pt]{amsart}
\usepackage{hyperref,fancyhdr,mathrsfs,amsmath,amscd,amsthm,amsfonts,latexsym,amssymb}

\DeclareMathOperator{\trace}{trace}

\DeclareMathOperator{\RicM}{^{\it M}\!Ric}
\DeclareMathOperator{\RicN}{^{\it N}\!Ric}
\DeclareMathOperator{\Ric}{Ric}
\DeclareMathOperator{\dif}{d}

\renewcommand{\H}{\mathscr{H}}
\newcommand{\V}{\mathscr{V}}
\newcommand{\F}{\mathscr{F}}
\newcommand{\Fa}{\mathcal{F}}
\newcommand{\tG}{\mathscr{G}}
\newcommand{\tD}{\mathscr{D}}

\DeclareMathOperator{\Bh}{{\it B}^{\H}}

\DeclareMathOperator{\Bv}{{\it B}^{\V}}

\def \a{\alpha}

\def \G{\Gamma}
\def \l{\lambda}
\def \o{\omega}

\def \phi{\varphi}
\def \p{\pi}
\def \r{\rho}

\def \R{\mathbb{R}}

\def \C{\mathbb{C}\,}

\def\widecheckg{g^{\hspace*{-2.5pt}\vbox to 5pt{\hbox to
0pt{\LARGE$\check{}$}}}\hspace*{2pt}}

\def\widecheckl{\lambda^{\hspace*{-3.5pt}\vbox to 8pt{\hbox to
0pt{\LARGE$\check{}$}}}\hspace*{2pt}}

\begin{document}

\title{Harmonic morphisms between Weyl spaces\\
and twistorial maps}
\author{E.~Loubeau and R.~Pantilie\;\dag}
\thanks{\dag\;Gratefully acknowledges that this work was partially supported by a
grant from the Conseil G\'en\'eral du Finist\`ere.}
\email{\href{mailto:Eric.Loubeau@univ-brest.fr}{Eric.Loubeau@univ-brest.fr},
       \href{mailto:Radu.Pantilie@imar.ro}{Radu.Pantilie@imar.ro}}
\address{E.~Loubeau, D\'epartement de Math\'ematiques, Laboratoire C.N.R.S. U.M.R. 6205,
Universit\'e de Bretagne Occidentale, 6, Avenue Victor Le Gorgeu, CS 93837,
29238 Brest Cedex 3, France}
\address{R.~Pantilie, Institutul de Matematic\u a ``Simion Stoilow'' al Academiei Rom\^ane,
Calea Grivi\c tei nr. 21, 010702, C.P. 1-764, 014700, Bucure\c sti, Rom\^ania}

\newtheorem{thm}{Theorem}[section]
\newtheorem{lem}[thm]{Lemma}
\newtheorem{cor}[thm]{Corollary}
\newtheorem{prop}[thm]{Proposition}

\theoremstyle{definition}

\newtheorem{defn}[thm]{Definition}
\newtheorem{rem}[thm]{Remark}
\newtheorem{exm}[thm]{Example}

\numberwithin{equation}{section}

\maketitle
\thispagestyle{empty}
\vspace{-0.5cm}
\section*{Abstract}
\begin{quote}
{\footnotesize  We show that Weyl spaces provide a natural context for
harmonic morphisms and we give the necessary and sufficient conditions under which
on an Einstein--Weyl space of dimension four there can be defined, locally, at least
five distinct foliations of dimension two which produce harmonic morphisms
(Theorem \ref{thm:five4to2}\,). Also, we describe the harmonic morphisms between
Einstein--Weyl spaces of dimensions four and three (Theorem \ref{thm:4to3}\,).}
\end{quote}

\section*{Introduction}

\indent
Harmonic morphisms between Riemannian manifolds are smooth maps which pull back (local)
harmonic functions to harmonic functions. By the basic characterisation theorem of
B.~Fuglede and T.~Ishihara, harmonic morphisms are harmonic maps which are horizontally
weakly conformal \cite{Fug}\,, \cite{Ish}\,.\\
\indent The simplest nontrivial examples of harmonic morphisms are
given by harmonic functions from a two-dimensional oriented
conformal manifold: any such harmonic morphism is, locally, the sum
of a $(+)$holomorphic function and a $(-)$holomorphic function (see
\cite{BaiWoo2}\,). Similar descriptions, in higher dimensions, can
be obtained if instead of $(\pm)$holomorphic functions we use the
more general notion of \emph{twistorial map} \cite{PanWoo-sd}\,. A
\emph{twistorial structure} on a complex manifold $M$ is given by a
foliation $\F$ on a complex manifold $P$ such that
$\F\cap{\rm ker}\dif\!\p=\{0\}$ where $\p:P\to M$ is a proper
complex analytic submersion; the leaf space of $\F$ is called the
\emph{twistor space} of the given twistorial structure $(P,M,\p,\F)$
and is usually denoted $Z(M)$\,.
It follows that, locally, we can find sections of $\p$ whose images
are foliated by leaves of $\F$\,; by projecting back through $\p$ we
endow $M$ with a sheaf of complex analytic submersions. Locally, a
\emph{twistorial function} is the composition of such a submersion
followed by a complex analytic function.\\
\indent
A sufficient condition for a complex analytic map $\phi:M\to N$ between
complex manifolds endowed with twistorial structures $(P,M,\p_P,\F)$ and
$(Q,N,\p_Q,\tG)$ to be twistorial, with respect to a map $\varPhi:P\to Q$\,,
is $\p_Q\circ\varPhi=\phi\circ\p_P$ and $\dif\!\varPhi(\F)=\tG$\,;
then, $\varPhi$ induces a map between the twistor spaces of
$(P,M,\p_P,\F)$ and $(Q,N,\p_Q,\tG)$ which we call the
\emph{twistorial representation} of $\phi$ (with respect to $\varPhi$).
To define the notion of twistorial map it is not necessary
to assume that $\F$ and $\tG$ are integrable; that is, we can speak
about twistorial maps between manifolds endowed with \emph{almost
twistorial structures}\, (see \cite[\S3]{PanWoo-sd} for the definition
of twistorial maps in full generality).\\
\indent
Thus, any twistorial structure on $M$ determines a sheaf of
twistorial functions $\Fa_M$ on $M$\,; such sheaves can be obtained
by complexifying the following examples of complex valued functions:\\
\indent
(i) $(\pm)$holomorphic functions on a Hermitian manifold; here $(M,c,J)$ is
the (germ-unique) complexification of a Hermitian manifold, $P=M_+\sqcup M_-$\,,
with $M_{\pm}=M$, and the twistor distribution $\F|_{M_{\pm}}$ is the
$(\mp{\rm i})$eigendistribution of $J$ (cf.\ \cite[Example 2.3]{PanWoo-sd}\,).\\
\indent
(ii) functions on a three-dimensional Einstein--Weyl space
which, locally, are the composition of a horizontally conformal submersion
with geodesic fibres followed by a complex analytic function; here
$(M^3,c,D)$ is the complexification of a
three-dimensional Einstein--Weyl space, $\p:P\to M$ is the bundle of
two-dimensional degenerate spaces on $(M^3,c)$ and, for each $p\in P$,
the subspace $\F_p\subseteq T_pP$ is the horizontal lift of
$p\subseteq T_{\p(p)}M$, with respect to $D$ (see \cite[Example 2.4]{PanWoo-sd}
for details about this twistorial structure).\\
\indent
(iii) functions on a four-dimensional anti-self-dual manifold which are holomorphic
with respect to a (local) positive Hermitian structure on it; here
$(M^4,c)$ is the complexification of a four-dimensional anti-self-dual
manifold, $\p:P\to M$ is the bundle of self-dual spaces on $(M^4,c)$
and, for each $p\in P$, the subspace $\F_p\subseteq T_pP$ is the horizontal
lift of $p\subseteq T_{\p(p)}M$, with respect to the Levi-Civita connection
of any local representative of $c$ (see \cite[Example 2.6]{PanWoo-sd} for
details about this twistorial structure).\\
\indent
One of the main steps in the process of classifying harmonic morphisms is to prove
that these are twistorial. For example, any harmonic morphism from an
Einstein manifold of dimension four with fibres of dimension one or two is twistorial
(\,\cite{Pan-4to3}\,, \cite{PanWoo-sd}\,, \cite{Woo-4d}\,; cf.\ Corollaries \ref{cor:4to2}
and \ref{cor:4to3}\,, below).\\
\indent
This paper gives an answer to the following question of John~C.~Wood:
\emph{Can (the submersive) twistorial maps be seen as harmonic morphisms}?
If we restrict ourselves to twistorial maps which pull back twistorial functions
to twistorial functions (for example, if we work with twistorial maps, as
above, with the further property that their twistorial representation is submersive)
then the answer,
in the affirmative, to this question follows if we work with sheaves of twistorial
functions $\Fa_M$ for which there exists a sheaf of `harmonic' functions $\mathcal{L}$
such that $\Fa_M\cap\mathcal{L}$ is a `sufficiently large' subsheaf of $\mathcal{L}$
(in particular, if the sheaf of vector spaces generated by $\Fa_M\cap\mathcal{L}$
is equal to $\mathcal{L}$\,, like in the case of two-dimensional conformal manifolds).
We argue that for each of the examples (i)\,, (ii)\,, (iii)\,, above, a good candidate for
$\mathcal{L}$ can be obtained by endowing the given conformal structure with
a suitable Weyl connection (the obvious one, for (ii)\,).\\
\indent
The definition of harmonic functions on a Weyl space
is given in Section \ref{section:harmorphWeyl}\,; there we also show that the
basic theorem of B.~Fuglede and T.~Ishihara
generalizes to harmonic morphisms between Weyl spaces.
In Section \ref{section:fundamn}\,, we do the same
for the fundamental equation for harmonic morphisms (see \cite{BaiWoo2}\,).\\
\indent
In Section \ref{section:Hermitian}\,,
we recall the definition and the basic properties of the
Weyl connection of an almost Hermitian manifold \cite{Vai}\,; we show that, for
Hermitian manifolds, this is characterized by the property that all the
$(\pm)$holomorphic functions are harmonic with respect to it. We also, prove
that any holomorphic horizontally weakly conformal map between almost Hermitian
manifolds endowed with their Weyl connections is harmonic and, hence,
a harmonic morphism (cf.\ \cite{Lic}\,, \cite{GudWoo}\,).\\
\indent
In Section \ref{section:hartwist}\,, we recall from \cite{PanWoo-sd}, the basic
examples of twistorial maps and we study the conditions under which these are
harmonic morphisms; the resulting harmonic morphisms are between Weyl spaces
of dimensions $m$ and $n$ where $(m,n)=(3,2)\,,\,(4,2)\,,\,(4,3)$\,. Thus,
we obtain the following for an analytic submersion $\phi$ between Weyl spaces of
dimensions $m$ and $n$\,.\\
\indent
$\bullet$ If $(m,n)=(3,2)$ then $\phi$ is twistorial if and only if it is a harmonic
morphism (Example \ref{exm:twistmaps3to2}\,).\\
\indent
$\bullet$ If $(m,n)=(4,2)\,,\,(4,3)$ then we obtain necessary and sufficient
conditions under which $\phi$ is a twistorial harmonic morphism
(Example \ref{exm:twistmaps4to2}\,, Theorem \ref{thm:twistharmorph4to3}\,).\\
\indent
The main purpose of Section \ref{section:Riccitwistsharmorphs} is to prove that
if $\phi$ is a harmonic morphism from a Weyl space of dimension four to a
Weyl space of dimension two or three then $\phi$ is
twistorial if and only if the difference between the Ricci tensor of the domain and
the pull-back of the Ricci tensor of the codomain is zero along the horizontal null
directions (Propositions \ref{prop:43to2} and \ref{prop:4to23}\,).\\
\indent
In Section \ref{section:4to23}\,, we apply the results of Sections \ref{section:hartwist}
and \ref{section:Riccitwistsharmorphs} to study harmonic morphisms from four-dimensional
Einstein--Weyl spaces. Thus,\\
\indent
$\bullet$ we give the necessary and sufficient conditions
under which on an Einstein--Weyl space of dimension four there can
be defined, locally, at least five distinct foliations of dimension
two which produce harmonic morphisms (Theorem \ref{thm:five4to2}\,);\\
\indent
$\bullet$ we describe the harmonic morphisms between
Einstein--Weyl spaces of dimensions four and three (Theorem \ref{thm:4to3}\,).\\

\indent
We are grateful to Paul~Baird for his interest in this work and to John~C.~Wood
for useful comments.

\section{Harmonic morphisms between Weyl spaces} \label{section:harmorphWeyl}

\indent
In this section we shall work in the smooth and (real or complex) analytic categories;
unless otherwise stated, all the manifolds are assumed to be connected.
A conformal manifold
$(M^m,c)$ is a manifold endowed with a reduction of its frame bundle to
$CO(m,\mathbb{K})$\,, ($\mathbb{K}=\R,\,\C$\,). We shall denote by $L^2$ the line
bundle associated to the bundle of conformal frames of $(M^m,c)$ via the morphism of
Lie groups $\rho_m:CO(m,\mathbb{K})\to\mathbb{K}\setminus\{0\}$ characterized by
$a^Ta=\rho_m(a)I_m$\,, ($a\in{}CO(m,\mathbb{K})\,,\,\mathbb{K}=\R,\,\C$\,).
In the smooth and real analytic categories $L^2$ is canonically oriented and so
it admits a canonical square root, denoted by $L$\,; moreover, $L$
does not depend of $c$ (see \cite{Cal-sds}\,). In the complex analytic category
such a square root can be found locally; furthermore, if $m$ is odd then there exists a
canonical choice for $L$ (this follows from the fact that, if $m$ is
odd then, there exists a natural isomorphism of Lie groups
$CO(m,\mathbb{K})=\mathbb{K}^*\times SO(m,\mathbb{K})$\,, ($\mathbb{K}=\R,\,\C$\,)\,).
In the smooth and real analytic categories positively oriented local sections
of $L^2$ correspond to local representatives of $c$\,. In the complex analytic
category, nowhere zero local sections of $L^2$ correspond to local representatives
of $c$\,. Note that, if $b$ is a section of $\otimes^2\,T^*\!M$ then its traces
with respect to local representatives of $c$ define a section of ${L^*}^{\,2}$
which will be denoted $\trace_cb$\,. More generally, if $b$ is a section of
$E\otimes(\otimes^2T^*M)$ for some vector
bundle $E$ over $M$ then we can define $\trace_cb$ which is a section of
$E\otimes {L^*}^{\,2}$\,;
if $E=TM$ then $(\trace_gb)^{\flat}$\,, where $g$ is any local
representative of $c$\,, defines a 1-form on $M$ which will be denoted
$(\trace_cb)^{\flat}$ (see \cite{Gau}\,, \cite{Cal-sds}\,).\\

\indent
The following definition is a simple generalization of the well-known notion of harmonic map (see \cite{BaiWoo2}\,).

\begin{defn} \label{defn:harwrt}
Let $(M,c)$ be a conformal manifold, $N$ a manifold and $D^M$, $D^N$ linear connections
on $M$, $N$, respectively.\\
\indent
A map $\phi:(M,c,D^M)\to(N,D^N)$ is called \emph{harmonic (with respect to $c$, $D^M$, $D^N$)} if
\begin{equation} \label{e:harwrt}
\trace_c(D\!\dif\!\phi)=0
\end{equation}
where $D$ is the connection on $\phi^*(TN)\otimes T^*M$ induced by $D^M$, $D^N$ and $\phi$\,.
\end{defn}

\indent
Obviously, there is no loss of generality if we assume $D^M$ and $D^N$ to be torsion free.\\
\indent
A harmonic map $\phi:(M,g)\to(N,h)$ between Riemannian manifolds is harmonic in the
sense of Definition \ref{defn:harwrt} if $M$ and $N$ are endowed with the Levi-Civita
connections of $g$ and $h$, respectively, and $M$ is considered with the conformal
structure determined by $g$\,.\\
\indent
We shall always consider $\mathbb{K}\,(=\R,\,\C)$ to be endowed with its conformal
structure and canonical connection (here $\C$ is considered to be a one-dimensional
complex manifold). Clearly, a curve on $(M,D)$ is harmonic if and only if it is a geodesic
of $D$\,.\\

\indent
Let $(M,c)$ be a conformal manifold. A torsion free conformal connection on
$(M,c)$ is called a
\emph{Weyl connection}; if $D$ is a Weyl connection on $(M,c)$ then $(M,c,D)$
is called a \emph{Weyl space}
(see \cite{Gau}\,). A function (locally) defined on a Weyl space $(M,c,D)$ will be called
\emph{harmonic} if it is harmonic with respect to $c$, $D$. If $\dim M=2$ then a function
$f$ on the Weyl space $(M,c,D)$
is harmonic if and only if it is harmonic with respect to any local representative of $c$\,.\\

\begin{prop}  \label{prop:samehar}
Let $(M,c_M)$ be a conformal manifold, of dimension $m\neq2$\,, endowed with a linear connection $D$\,.\\
\indent
Then there exists a unique Weyl connection $D_1$ on $(M,c_M)$ such that
\begin{equation} \label{e:samehar}
\trace_{c_M}(D\!\dif\!f)=\trace_{c_M}(D_1\!\dif\!f)
\end{equation}
for any function $f$ (locally) defined on $M$.
\end{prop}
\begin{proof}
For each local representative $g$ of $c_M$ we define a (local) 1-form $\a^g$ by
\begin{equation} \label{e:sameharLee}
\a^g(X)=\frac{1}{m-2}\,g(\trace_g(\nabla^g-D),X)
\end{equation}
for all $X\in TM$, where $\nabla^g$ is the Levi-Civita connection
of $g$\,. It is easy to prove that $\a^{g\l^{-2}}=\a^g+\l^{-1}\!\dif\!\l\,.$
Hence, the family of 1-forms $\{\a^g\}$ defines a connection
on $L$\,. But any connection on $L$ corresponds to a Weyl connection $D_1$ on $(M,c_M)$ (see \cite{Gau}\,, the 1-form
$\a^g$ is \emph{the Lee form of $D_1$ with respect to $g$}\,).
Now, \eqref{e:samehar} is equivalent to \eqref{e:sameharLee} and the proof follows.
\end{proof}

\indent
The following definition (cf.\ \cite{Fug}\,, \cite{Ish}\,, \cite{BaiWoo2}\,) will be
central in this paper.

\begin{defn} \label{defn:harmorphw}
Let $(M,c_M,D^M)$ and $(N,c_N,D^N)$ be Weyl manifolds.\\
\indent A map $\phi:(M,c_M,D^M)\to(N,c_N,D^N)$ is called a
\emph{harmonic morphism} if for any harmonic function $f$ defined on
some open set $U$ of $N$, such that $\phi^{-1}(U)\neq\emptyset$\,,
the function $f\circ\phi|_{\phi^{-1}(U)}$ is harmonic.
\end{defn}

\begin{rem}
Proposition \ref{prop:samehar} shows that, if $\dim M,\dim N\neq2$ then
Definition \ref{defn:harmorphw} does not become more general by using linear
connections instead of Weyl connections.
\end{rem}

\indent
Any harmonic morphism between Riemannian manifolds $\phi:(M,g)\to(N,h)$ is also
a harmonic morphism between Weyl spaces $\phi:(M,[g],\nabla^g)\to(N,[h],\nabla^h)$\,,
where $[g]$\,, $[h]$
are the conformal structures determined by $g$\,, $h$ and $\nabla^g$\,, $\nabla^h$
are the Levi-Civita connections
of $g$\,, $h$\,, respectively. However, not all harmonic morphisms between Weyl
spaces arise in this way (see, for example, Remark \ref{rem:fundamn}(2)\,,
Corollary \ref{cor:holharmorph} and Example \ref{exm:GibHaw,Bel,Hit}(2)\,, below).\\

\begin{defn}[cf.\ \cite{BaiWoo2}\,]
Let $(M,c_M,D^M)$ be a Weyl space. A foliation $\V$ on $M$ \emph{produces harmonic
morphisms} on $(M,c_M,D^M)$ if, locally, it can be defined by submersive harmonic morphisms;
that is, for any point of $x\in M$ there exists a submersive harmonic
morphism $\phi:(U,c_M|_U,D^M|_U)\to(N,c_N,D^N)$ such that
$\V|_U={\rm ker}\dif\!\phi$\,, where $U$ is an open neighbourhood of $x$\,.
\end{defn}

\indent
Next we shall prove the Fuglede-Ishihara theorem (\cite{Fug}\,, \cite{Ish}\,, see
\cite{BaiWoo2}\,) for harmonic morphisms
between Weyl spaces. For this we apply the standard strategy (see \cite{BaiWoo2}\,)\,.
Firstly, we need an existence result for harmonic functions from Weyl spaces.

\begin{lem}[cf.\ \cite{BaiWoo2}\,]  \label{lem:harexist}
Let $(M,c,D)$ be a Weyl space and let $x\in M$.\\
\indent
Then for any $v\in T^*_xM$ and any trace free symmetric bilinear form $b$ on $(T_xM,c_x)$
there exists a harmonic function $f$ defined on some open neighbourhood of $x$
such that $\dif\!f_x=v$ and $(D\!\dif\!f)_x=b$\,.
\end{lem}
\begin{proof}
This is essentially the same as for harmonic functions on Riemannian manifolds
(see \cite{BaiWoo2} and the references therein).\\
\indent
We shall give a straightforward proof assuming $(M,c,D)$ (real or complex) analytic
(cf.\ \cite[Lemma 2]{Lou}\,, where we
assume that the metric is analytic).
Let $U$ be the domain of a normal coordinate system $x^1,\ldots,x^m$ for $D$, centred
at $x$\,, where $m=\dim M$.
We may assume $g(\dif\!x^m,\dif\!x^m)=1$\,, at $x$, for some local
representative $g$ of $c$ over $U$. Hence, by passing to a smaller open neighbourhood
of $x$\,, if necessary, we may
assume that the hypersurface $S=\{x^m=0\}$ is nowhere degenerate; equivalently, $S$
is noncharacteristic
for the second order differential operator $f\mapsto\trace_g(D\!\dif\!f)$\,.\\
\indent
Let $p=b_{ij}x^ix^j+v_ix^i$\,. Then, by further restricting $U$, if necessary,
and by applying the Cauchy-Kovalevskaya theorem,
we can find a harmonic function $f$, with respect to $c,D$, defined on $U$ such
that $f$ and $p$ are equal up to the
first derivatives along $S$\,; in particular, $\dif\!f_x=v$\,. Hence, possibly
excepting $\frac{\partial^2\!f}{(\partial x^m)^2}(x)$\,,
all the second order partial derivatives of $f$, at $x$\,, are equal to the
corresponding derivatives of $p$\,, at $x$\,.
As $f$ is harmonic, $b$ is trace free, with respect to $g$\,, and $x$ is the
centre of the normal system of coordinates
$x^1,\ldots,x^m$\,, for $D$\,, the derivatives
$\frac{\partial^2f}{(\partial x^m)^2}(x)$ and $\frac{\partial^2p}{(\partial x^m)^2}(x)$
are determined by the other second order partial derivatives, at $x$\,, of $f$ and $p$\,,
respectively, and hence
must be equal. Thus $(D\!\dif\!f)_x=b$\,.
\end{proof}

\begin{rem}
Let $f$ be a harmonic function (locally defined) on a Weyl space $(M,c,D)$ and let
$x\in M$ such that $\dif\!f_x\neq0$\,.
Then there exists a local representative $g$ of $c$ defined on some neighbourhood
$U$ of $x$ such that $f$ is harmonic with
respect to $g$ (this follows, for example, from \eqref{e:sameharLee} applied to $D$\,).
However, the sheaf of harmonic
functions on $U$\,, with respect to $c$, $D$\,, is equal to the sheaf of harmonic
functions of $g$ if and only if
$D$ is the Levi-Civita connection of $g$\,.
\end{rem}

\indent
Let $(M,c_M)$ be a conformal manifold and let $L_M$ be the associated line bundle.
Note that, the conformal structure $c_M$ corresponds to a `musical' isomorphism
$T^*\!M\stackrel{\sim}{\longrightarrow}TM\otimes L_M^*\!^2$\,.
Hence, the differential of any map $\phi:(M,c_M)\to(N,c_N)$ between conformal
manifolds has an adjoint $(\dif\!\phi)^T$ which is a section of the vector bundle
${\rm Hom}\bigl(\phi^*(TN),TM\bigr)\otimes
{\rm Hom}\bigl(L_M^{\,2},\phi^*(L_N^{\,2}\,)\bigr)$\,.

\begin{defn}[cf.\ \cite{BaiWoo2}\,, \cite{Pamb}\,] \label{defn:confhorconf}
A map $\phi:(M,c_M)\to(N,c_N)$ between conformal manifolds is called
\emph{horizontally weakly conformal} if there exists a section
$\Lambda$ of ${\rm Hom}\bigl(L_M^{\,2},\phi^*(L_N^{\,2}\,)\bigr)$ such that
$$((\dif\!\phi)^T)^*(c_M)=\Lambda\otimes c_{\phi^*(TN)}$$
where $c_{\phi^*(TN)}$ is the conformal structure on $\phi^*(TN)$
induced by $\phi$ and $c_N$\,; $\Lambda$ is called the
\emph{square dilation} of $\phi$\,. If $\Lambda=\l^2$ for some section $\l$ of
${\rm Hom}(L_M,\phi^*(L_N))$ then
$\l$ is a \emph{dilation} of $\phi$\,.\\
\indent
Let $(M,c_M)$ be a conformal manifold. A \emph{conformal foliation} on $(M,c_M)$
is a foliation which is locally defined by horizontally conformal submersions.
\end{defn}

\indent
In the complex analytic category, let $\phi:(M,c_M)\to(N,c_N)$ be a horizontally conformal
submersion with nowhere degenerate fibres. Then each point point of $M$ has an open
neighbourhood $U$ such that $\phi|_U$ admits a dilation. In fact, by passing, if necessary,
to a conformal $\mathbb{Z}_2$-covering space of $M$ we can suppose that $\phi$ admits a
dilation.\\
\indent
In the smooth (real analytic) category any horizontally
conformal submersion admits a dilation. Furthermore, if
$\phi:(M,c_M)\to(N,c_N)$ is horizontally weakly conformal then there
exists a continuous section $\l$ of ${\rm Hom}(L_M,\phi^*(L_N))$
such that $\Lambda=\l^2$\,; moreover, $\l$ is smooth (real analytic)
over the set of regular points of $\phi$ (cf.\ \cite{BaiWoo2}\,).\\

\indent
Now we can prove the Fuglede-Ishihara theorem for harmonic morphisms between Weyl spaces.

\begin{thm}[cf.\ \cite{Fug}\,, \cite{Ish}\,]  \label{thm:Fug-Ish}
Let $(M,c_M,D^M)$ and $(N,c_N,D^N)$ be Weyl manifolds.\\
\indent
A map $\phi:(M,c_M,D^M)\to(N,c_N,D^N)$ is a harmonic morphism if and only if it is harmonic, with respect to $c_M$, $D^M$, $D^N$, and horizontally weakly conformal.
\end{thm}
\begin{proof}
For any function $f$, locally defined on $M$, a straightforward calculation gives
\begin{equation} \label{e:chain}
\trace_{c_M}\bigl(D\!\dif(f\circ\phi)\bigr)=\dif\!f\bigl(\trace_{c_M}(D\!\dif\!\phi)\bigr)+c_N\bigl(D\!\dif\!f,((\dif\!\phi)^T)^*(c_M)\bigr)\;.
\end{equation}
\indent
By applying Lemma \ref{lem:harexist} with $b=0$ and for all $v\in T_x^*N$, ($x\in N$), from equation \eqref{e:chain} we obtain that $\phi$
is harmonic with respect to $c_M$, $D^M$, $D^N$. Then, by applying again
Lemma \ref{lem:harexist}\,, from equation \eqref{e:chain}
we obtain that for all trace free symmetric $b\in\otimes^2T^*N$ we have
$c_N\bigl(b,((\dif\!\phi)^T)^*(c_M)\bigr)=0$\,. The proof follows.
\end{proof}

\indent
An immediate consequence of Theorem \ref{thm:Fug-Ish} is that any foliation which produces
harmonic morphisms is a conformal foliation (cf.\ \cite{BaiWoo2}\,).

\section{The fundamental equation}  \label{section:fundamn}

\indent
In this section we shall work in the smooth and (real or complex) analytic categories.
Let $(M^m,c_M)$ be a conformal manifold and let $L$ be the associated line bundle on $M$.
Let $\V\subseteq TM$ be a nondegenerate distribution, and let $\H=\V^{\perp}$ be its
orthogonal complement.  Then $c_M$ induces conformal structures $c_M|_{\V}$ and $c_M|_{\H}$
on $\V$ and $\H$, respectively. Let $L_{\V}$ and $L_{\H}$ be the line bundles on $M$
determined by the conformal structures $c_M|_{\V}$ and $c_M|_{\H}$\,, respectively.
As any local representative of $c_M$ induces local representatives of the conformal
structures induced on $\V$ and $\H$, we have isomorphisms
between $L^2$\,, $L_{\V}^2$ and $L_{\H}^2$ (seen as bundles with group
$\bigl((0,\infty),\cdot\bigr)$\,, in the smooth and real analytic categories);
we shall always identify $L^2=L_{\V}^2=L_{\H}^2$\,, in this way. Conversely, conformal
structures on the complementary distributions $\V$ and $\H$ together with an isomorphism
between $L_{\V}^2$ and $L_{\H}^2$
determine a conformal structure on $M$ such that $\H=\V^{\perp}$ \cite{Cal-sds}\,.
In other words, nondegenerate distributions $\V$\,, of dimension $m-n$\,, on
$(M^m,c_M)$ correspond to reductions of $c_M$ to the subgroup
$G=\bigl\{\,(a,b)\in CO(m-n)\times CO(n)\,|\,\r_{m-n}(a)=\r_n(b)\,\bigr\}$ of $CO(m)$\,.
Then, as the morphisms of Lie groups $p_1:G\to CO(m-n)$\,, $(a,b)\mapsto a$\,,
and $p_2:G\to CO(n)$\,, $(a,b)\mapsto b$\,, satisfy
$\r_m|_G=\r_{m-n}\circ p_1=\r_n\circ p_2$\,, we obtain that
$c_M|_{\V}=p_1(c_M)$ and $c_M|_{\H}=p_2(c_M)$ are such that $L_{\V}^2$ and $L_{\H}^2$ are
isomorphic to $L^2$\,. Conversely, if $c_{\V}$ and $c_{\H}$ are (the bundles of conformal
frames of) conformal structures
on the complementary distributions $\V$ and $\H$, respectively, then $c_{\V}+c_{\H}$
is a reduction of the bundle of linear frames on $M^m$ to $CO(m-n)\times CO(n)$ and,
it is easy to see that, isomorphisms
between $L_{\V}^2$ and $L_{\H}^2$ correspond to reductions of $L_{\V}^2\oplus L_{\H}^2$
to $\iota:H\hookrightarrow H\times H$\,, $a\mapsto(a,a)$\,, where
$H=\bigl((0,\infty),\cdot\bigr)$ in the smooth and real analytic categories,
and $H=\bigl(\C\setminus\{0\},\cdot\bigr)$ in the complex analytic category.
As $G=(\r_{m-n}\times\r_n)^{-1}(\iota(H))$\,,
it follows that reductions of $c_{\V}+c_{\H}$ to $G$ correspond to isomorphisms between
$L_{\V}^2$ and $L_{\H}^2$\,; any such reduction determines a conformal structure $c_M$
on $M^m$ such that $\H=\V^{\perp}$ and $c_M|_{\V}=c_{\V}$\,, $c_M|_{\H}=c_{\H}$\,.

\begin{exm}[\,\cite{Bott-partial}\,] \label{exm:Bott-partial}
\indent
Let $M$ be a manifold endowed with two complementary distributions $\V$ and $\H$.
The \emph{Bott partial connection} $D^{\rm Bott}$ on $\V$\,,
over $\H$, is defined by $D^{\rm Bott}_XU=\V[X,U]$ for local sections $X$ of $\H$ and
$U$ of $\V$.\\
\indent
Suppose that $M$ is endowed with a conformal structure $c_M$ with respect to which
$\V$ is nondegenerate
and $\H=\V^{\perp}$. As $(L^2)^{m-n}=(\Lambda^{m-n}\V)^2$\,, where $n$ is the dimension
of the distribution $\H$,
$D^{\rm Bott}$ induces a partial connection on $L$ which will also be denoted
$D^{\rm Bott}$\,; the local connection form of this connection with respect to a
local section of $L$\,, corresponding
to a local representative $g$ of $c_M$\,, is $\frac{1}{m-n}\trace_g(\Bv)^{\flat}$.
\end{exm}

\indent
Let $\phi:(M,c_M,D^M)\to(N,c_N,D^N)$ be a horizontally conformal submersion with nowhere degenerate
fibres between Weyl spaces. We shall denote, as usual (see, for example, \cite{BaiWoo2}\,),
$\V={\rm ker}\,\dif\!\phi$\,, $\H=\V^{\perp}$\,. Then $D^M$ and $D^N$ induce Weyl partial connections,
with respect to $\V$, on $(\H,c_M|_{\H})$\,, over $\H$, which will be denoted $\H D^M$ and $D^N$\,,
respectively. (Recall (see \cite{PanWoo-sd}\,) that a Weyl partial connection $D$ on $(\H,c)$\,,
over $\H$, is a conformal partial connection $D$ on $(\H,c)$ whose torsion tensor field $T$,
with respect to $\V$, defined by $T(X,Y)=D_XY-D_YX-\H[X,Y]$ for local sections $X$ and $Y$ of $\H$, is zero.)\\
\indent
If $D$ is a (partial) connection on $L$ and $k\in\mathbb{Z}$ then we shall denote by $D^k$ the
(partial) connection induced on $L^k$ where, $L^k=\otimes^kL$ if $k$ is a natural number and  $L^k=\otimes^{-k\,}L^*$ if $k$ is a negative integer.

\begin{prop} \label{prop:fundamn}
Let $\phi:(M^m,c_M,D^M)\to(N^n,c_N,D^N)$ be a horizontally conformal submersion with nowhere degenerate
fibres between Weyl spaces. Then
\begin{equation} \label{e:fundamn}
\trace_{c_M}(D\!\dif\!\phi)^{\flat}=(\H D^M)^{m-2}\otimes(D^N)^{-(n-2)}-(D^{\rm Bott})^{m-n}\;.
\end{equation}
\end{prop}
\begin{proof}
Let $B^{\V,D^M}$ be the second fundamental form of $\V$, with respect to $D^M$, defined by
$$B^{\V,D^M}(U,V)=\frac12\,\V(D_UV+D_VU)$$ for local sections $U$ and $V$ of $\V$ (see \cite{Cal-sds}\,,
cf.\ \cite{BaiWoo2}\,). A straightforward calculation gives
\begin{equation} \label{e:fundamn1}
\trace_{c_M}(D\!\dif\!\phi)^{\flat}=\trace_{c_M}(D^N-\H D^M)-\trace_{c_M}(B^{\V,D^M})^{\flat}\;.
\end{equation}
\indent
Now let $g$ be a local representative of $c_M$\,, corresponding to some local section $s$ of $L$\,,
and let $\a_M$ and $\a_N$ be the Lee forms
of $D^M$ and $D^N$, respectively, with respect to $g$\,. Recall (see \cite{Gau}\,) that
$\a_M$ ($\a_N$) is the local connection form of $D^M$ ($D^N$) with respect to $s$\,.
Also, it is easy to prove that
\begin{equation} \label{e:traceBVD}
\trace_g(B^{\V,D^M})^{\flat}=\trace_g(\Bv)^{\flat}-(m-n)\a_M|_{\H}
\end{equation}
where $\Bv$ is the second fundamental form of $\V$ with respect to
(the Levi-Civita connection of) $g$\,.\\
\indent
It follows that \eqref{e:fundamn1} is equivalent to
\begin{equation} \label{e2fundamn2}
\trace_{c_M}(D\!\dif\!\phi)^{\flat}=(m-2)\a_M|_{\H}-(n-2)\a_N-\trace_g(\Bv)^{\flat}
\end{equation}
\indent
The proof follows from Example \ref{exm:Bott-partial}\,.
\end{proof}

\begin{rem} \label{rem:fundamn}
1) When $D^M$ and $D^N$ are the Levi-Civita connections of (local) representatives
of $c_M$ and $c_N$\,, respectively, then \eqref{e:fundamn} reduces to the
\emph{fundamental equation} (see \cite{BaiWoo2}\,) for horizontally conformal
submersions between Riemannian manifolds.\\
\indent
2) Let $\phi:(M,c_M)\to(N,c_N,D^N)$ be a horizontally conformal submersion,
with nowhere degenerate fibres, from a conformal manifold to a Weyl space.
{}From the fundamental equation \eqref{e:fundamn}\,, it follows that there exists
a Weyl connection $D^M$ on $(M,c_M)$ such that $\phi:(M,c_M,D^M)\to(N,c_N,D^N)$ is
a harmonic morphism.\\
\indent
3) Another consequence of Proposition \ref{prop:fundamn} is that if
$D_1$ and $D_2$ are Weyl connections on a conformal manifold $(N,c_N)$\,, of dimension
not equal to two, and $\phi:(M,c_M,D^M)\to(N,c_N,D_j)$\,, $j=1\,,\,2$\,, is a surjective
harmonic morphism with nowhere degenerate fibres then $D_1=D_2$ (cf.\ \cite{BaiWoo2}\,).
\end{rem}

\indent
We shall say that $\V$ is \emph{minimal, with respect to $D^M$,} if
$\trace_{c_M}(B^{\V,D^M})=0$\,; then $\V$
is minimal, with respect to $D^M$, if and only if $D^M$ and $D^{\rm Bott}$
induce the same partial connection on $L$\,, over $\H$ (cf.\ \cite{Cal-sds}\,).\\
\indent
Similar to the case of harmonic morphisms between Riemannian manifolds, from the
\emph{fundamental equation} \eqref{e:fundamn}\,, we obtain the following.

\begin{thm}[cf.\ \cite{BaiEel}\,]  \label{thm:BaiEel}
Let $\phi:(M,c_M,D^M)\to(N,c_N,D^N)$ be a horizontally conformal submersion with
nowhere degenerate fibres between Weyl spaces.\\
\indent
{\rm (a)} If $\dim N=2$ then $\phi$ is a harmonic morphism if and only if its
fibres are minimal, with respect to $D^M$.\\
\indent
{\rm (b)} If $\dim N\neq2$ then any two of the following assertions imply the third:\\
\indent\quad{\rm (i)} $\phi$ is a harmonic morphism.\\
\indent\quad{\rm (ii)} The fibres of $\phi$ are minimal, with respect to $D^M$.\\
\indent\quad{\rm (iii)} $\H D^M=D^N$\,.
\qed
\end{thm}

\indent
Let $(M,c_M)$ be a conformal manifold endowed with a
conformal foliation $\V$. Note that, if $\V$ is nowhere degenerate
then $L$ is basic (with respect to $\V$\,). Indeed, for any point
$x\in M$ there exists an open neighbourhood $U$, of $x$\,, on which
there can be defined a horizontally conformal submersion
$\phi:(U,c_M|_U)\to(N,c_N)$ which defines $\V$ on $U$ and which
admits a dilation $\l$\,. Then $\l$ is an isomorphism between
$L|_U$ and $\phi^*(L_N)$.\\
\indent
The following proposition follows easily from \eqref{e:fundamn}\,.

\begin{prop}[cf.\ \cite{BaiWoo2}\,]
Let $(M^m,c_M,D^M)$ be a Weyl space, $\dim M=m$\,, endowed with a nowhere
degenerate conformal foliation $\V$.\\
\indent
{\rm (i)} If ${\rm codim}\V=2$ then $\V$ produces harmonic morphisms on $(M,c_M,D^M)$
if and only if its leaves are minimal.\\
\indent
{\rm (ii)} If ${\rm codim}\V=n\neq2$ then $\V$ produces harmonic morphisms on $(M,c_M,D^M)$
if and only if the partial connection $(\H D^M)^{m-2}\otimes(D^{\rm Bott})^{-(m-n)}$
on $L^{n-2}$\,, over $\H$, is basic.
\end{prop}

\indent
We end this section with an example of a Weyl connection which will be useful later on.

\begin{exm}[\,\cite{Cal-sds}\,] \label{exm:minimalWeyl}
Let $(M^m,c)$ be a conformal manifold endowed with a nondegenerate distribution $\V$,
of codimension $n$\,, and let $\H=\V^{\perp}$.\\
\indent
For each local representative $g$ of $c_M$ define a (local) 1-form $\a^g$ by
\begin{equation} \label{e:minimalWeyl}
\a^g=\frac{1}{m-n}\,\trace_g(\Bv)^{\flat}+\frac{1}{n}\,\trace_g(\Bh)^{\flat}\;.
\end{equation}
Then $\a^{g\l^{-2}}=\a^g+\l^{-1}\!\dif\!\l$\,. Hence, the family of 1-forms
$\{\a^g\}$ defines a Weyl connection $D$ on $(M^m,c)$\,. The Weyl connection $D$
is called the \emph{(minimal) Weyl connection of $(M^m,c,\V)$}\,. Note that,
if we denote by $D^{{\rm Bott},\V}$ and $D^{{\rm Bott},\H}$ the partial connections,
over $\H$ and $\V$, respectively, induced on $L$ by the Bott partial connections
(see Example \ref{exm:Bott-partial}\,) on $\V$ and $\H$, then
$D=D^{{\rm Bott},\V}+D^{{\rm Bott},\H}$. Also, note that
$\trace_c(D\V)=\trace_c(D\H)=0$ (cf.\ Remark \ref{rem:PWeyl}\,, below).\\
\indent
As \eqref{e:traceBVD} holds without the assumption that $\V$ is conformal, $\V$
and $\H$ are minimal with respect to $D$\,; moreover, $D$ is the unique Weyl connection
on $(M^m,c)$ with this property. It follows that if $\V$ is one-dimensional and
conformal then the connection induced by $D$ on $L$ is flat if and only if $\V$ is
locally generated by conformal vector fields.
\end{exm}

\section{Harmonic maps and morphisms between almost Hermitian manifolds} \label{section:Hermitian}

\indent
In this section we shall work in the smooth and (real or complex) analytic categories.
An \emph{almost Hermitian (conformal) manifold} is a triple $(M,c,J)$ where
$(M,c)$ is a conformal
manifold and $J$ is a compatible almost complex structure; that is, if we consider $c$
as an $L^2$-valued Riemannian metric on $M$ \cite{Gau} then we have $c(JX,JY)=c(X,Y)$\,,
$(X,Y\in TM)$\,. Therefore, $\dim M$ is even and the \emph{K\"ahler form of $(M,c,J)$}\,,
defined by, $\o(X,Y)=c(JX,Y)$\,, $(X,Y\in TM)$\,, is an $L^2$-valued almost
symplectic structure
on $M$. A \emph{Hermitian (conformal) manifold} is an almost Hermitian
manifold $(M,c,J)$ such that $J$ is integrable.\\
\indent
To any almost Hermitian manifold, of dimension at least four, can be associated,
in a natural way, a Weyl connection, as follows.

\begin{prop}[\,\cite{Vai}\,] \label{prop:canonicalWeyl}
Let $(M,c,J)$ be an almost Hermitian manifold, of dimension $m\geq4$\,, and let
$\o\in\G(L^2\otimes\Lambda^2T^*M)$ be its K\"ahler form.\\
\indent
There exists a unique Weyl connection $D$ on $(M,c)$ such that $\trace_c(DJ)=0$\,,
the Lee form of $D$ with respect to a local representative $g$ of $c$\,,
is equal to $-\frac{1}{m-2}$ times the Lee form of $J$ with respect to $g$\,.
\end{prop}
\begin{proof}
Let $m=2n$\,, $(n\geq2)$\,. {}From the fact that $\o$ is an $L^2$-valued almost symplectic
structure on $M$, it follows (see \cite{Cal-sds}\,) that there exists a unique connection
$D$ on $L^2$ such that
\begin{equation} \label{e:canonicalWeyl}
\dif^D\!\o\wedge\o^{n-2}=0\;.
\end{equation}
We shall denote by the same letter $D$ the induced connection on $L$ and the corresponding
Weyl connection on $(M,c)$\,. Let $s$ be a local section of $L$ and let $\o^s$ be the
K\"ahler form of $J$ with respect to the local representative $g^s$ of $c$ corresponding
to $s$\,; that is, $\o^s(X,Y)=g^s(JX,Y)$\,, $(X,Y\in TM)$\,. It is easy to prove that
\eqref{e:canonicalWeyl} is equivalent to the fact that, for any local section
$s$ of $L$\,, the local connection form of $D$\,, with respect to $s$\,, is equal
to $-\frac{1}{m-2}$ times the Lee form of $J$\,, with respect to $g^s$\,.\\
\indent
Furthermore, \eqref{e:canonicalWeyl} is also equivalent to
$\sum_{i=1}^{n}(\dif^D\!\o)(X_i,JX_i,\cdot)=0$ for any conformal frame
$\bigl\{X_1,JX_1,\dots,X_n,JX_n\bigr\}$\,. Therefore to end the proof it is sufficient
to show that for any Weyl connection $D$ on $(M,c)$ we have
$$\sum_{i=1}^{n}(\dif^D\!\o)(X_i,JX_i,JY)=-c(\trace_g(DJ),Y)$$
for any $Y\in TM$ and where $g$ is the metric determined by the conformal frame
$\bigl\{X_1,JX_1,\cdots,X_n,JX_n\bigr\}$\,.
Indeed, as $Dc=0$ we have $(D\o)(X,Y)=c((DJ)(X),Y)$ and $(D\o)(X,JX)=0$\,,
$(X,Y\in TM)$\,. Therefore
\begin{equation*}
\begin{split}
\sum_{i=1}^{n}(\dif^D\!\o)(X_i,JX_i,JY)&=(D_{X_i}\o)(JX_i,JY)+(D_{JX_i}\o)(JY,X_i)+(D_{JY}\o)(X_i,JX_i)\\
&=c((D_{X_i}J)(JX_i),JY)+c((D_{JX_i}J)(JY),X_i)\\
&=-c(J(D_{X_i}J)(X_i),JY)-c(J(D_{JX_i}J)(Y),X_i)\\
&=-c((D_{X_i}J)(X_i),Y)-c(Y,(D_{JX_i}J)(JX_i))\\
&=-c(\trace_g(DJ),Y)\;.
\end{split}
\end{equation*}
\end{proof}

\begin{defn}[\,\cite{Vai}\,] \label{defn:canonicalWeyl}
Let $(M,c,J)$ be an almost Hermitian manifold, of dimension $\dim M\geq4$\,.\\
\indent
The \emph{Weyl connection of $(M,c,J)$} is the Weyl connection $D$ on $(M,c)$
such that $\trace_c(DJ)=0$\,.
\end{defn}

\begin{rem} \label{rem:canonicalWeyl}
1) \cite{Vai} Let $(M,c,J)$ be an almost Hermitian manifold, of dimension
$\dim M\geq4$\,, and let
$D$ be a Weyl connection on $(M,c)$\,.\\
\indent
Let $\nabla$ be the Levi-Civita connection of a local representative $g$ of $c$\,.
Then $D_{JX}J-J\,D_XJ=\nabla_{JX}J-J\,\nabla_XJ$\,, $(X\in TM)$\,. Hence $J$ is integrable
if and only if $D_{JX}J=J\,D_XJ$\,, $(X\in TM)$\,.\\
\indent
On the other hand, the condition $D_{JX}J=-J\,D_XJ$\,, $(X\in TM)$\,, is equivalent to
$(\dif^D\!\o)^{(1,2)\oplus(2,1)}=0$ and is a sufficient condition for $D$ to be the
Weyl connection of $(M,c,J)$\,. Hence, if $\dim M=4$ then $D_{JX}J=-J\,D_XJ$\,, $(X\in TM)$\,,
if and only if $D$ is the Weyl connection of $(M,c,J)$\,.\\
\indent
Thus, if $\dim M=4$ then $DJ=0$ if and only if $J$ is integrable and $D$ is the Weyl
connection of $(M,c,J)$\,. If $\dim M\geq6$ then it follows that $DJ=0$ if and only if,
locally, there exist representatives $g$ of $c$ with respect to which $(M,g,J)$ is
K\"ahler.\\
\indent
2) Let $(M,c,J)$ be an almost Hermitian manifold and let $f$ be a $(\pm)$holomorphic
function locally defined on $(M,J)$\,. (If $(M,c,J)$ is complex analytic then by a
$(\pm)$holomorphic function we mean a function which is constant along curves tangent
to the $(\mp{\rm i})$eigendistributions of $J$.)
If $\dim M=2$ then $f$ is harmonic with respect to any local representative
of $c$ (see \cite{BaiWoo2}\,). If $\dim M\geq4$ and $D$ is the Weyl connection of $(M,c,J)$
then $f$ is a harmonic function of $(M,c,D)$\,.\\
\indent
Furthermore, if $(M,c,J)$ is a Hermitian manifold, $\dim M\geq4$\,, then for any
Weyl connection $D$ on $(M,c)$ the following assertions are equivalent:\\
\indent
\quad(i) $D$ is the Weyl connection of $(M,c,J)$\,.\\
\indent
\quad(ii) Any $(\pm)$holomorphic function of $(M,J)$ is a harmonic function of $(M,c,D)$\,.\\
\indent
See Proposition \ref{prop:canonicalWeyl2} for a reformulation of this equivalence,
in the complex analytic category.
\end{rem}

\indent
Next, we prove the following useful lemma.

\begin{lem}[cf.\ \cite{Lic}\,, \cite{PanWoo-d}\,] \label{lem:A3}
Let $D^M$, $D^N$ be torsion free connections on the almost complex manifolds $(M,J^M)$,
$(N,J^N)$, respectively. Suppose that $M$ is endowed with a conformal structure $c$
and let $\phi:(M,J^M)\to(N,J^N)$ be a holomorphic map. Then
\begin{equation} \label{e:A3}
\trace_c\phi^*(D^N\!J^N)-\dif\!\phi(\trace_c(D^M\!J^M))+J^N(\trace_c(D\!\dif\!\phi))=0\;.
\end{equation}
\end{lem}
\begin{proof}
It is easy to prove that, for $X,Y\in TM$, we have
\begin{equation*}
D\!\dif\!\phi(X,J^MY)=(D^N_{\,\dif\!\phi(X)}J^N)(\dif\!\phi(Y))-\dif\!\phi((D^M_{\;X}J^M)(Y))+J^N(D\!\dif\!\phi(X,Y))\;.
\end{equation*}
The proof follows.
\end{proof}

{}From Lemma \ref{lem:A3} we easily obtain the following proposition
(cf.\ Remark \ref{rem:canonicalWeyl}(2)\,).

\begin{prop}[cf.\ \cite{Lic}\,, \cite{GudWoo}\,] \label{prop:Lic}
Let $(M,c_M,J^M)$ and $(N,c_N,J^N)$ be almost Hermitian manifolds. If
$\dim M\geq4$\,, $\dim N\geq4$ let $D^M$\,, $D^N$ be the Weyl connections
of $(M,c_M,J^M)$\,, $(N,c_N,J^N)$\,, respectively; if $\dim M=2$
or $\dim N=2$ then $D^M$ or $D^N$ will denote any Weyl connection on
$(M,c_M)$ or $(N,c_N)$\,, respectively.\\
\indent
Let $\phi:(M,J^M)\to(N,J^N)$ be a holomorphic map.\\
\indent
\quad{\rm (i)} If $(\dif^{D^N}\!\o_N)^{(1,2)\oplus(2,1)}=0$\,, where $\o_N$ is
the K\"ahler form of $(N,c_N,J^N)$\,, then $\phi:(M,c_M,D^M)\to(N,c_N,D^N)$
is a harmonic map.\\
\indent
\quad{\rm (ii)} If the map $\phi:(M,c_M)\to(N,c_N)$ is horizontally weakly conformal
and has nowhere degenerate fibres then $\phi:(M,c_M,D^M)\to(N,c_N,D^N)$
is a harmonic map and hence a harmonic morphism. \qed
\end{prop}

\indent
Note that, in assertion (i) of Proposition \ref{prop:Lic}\,, if $\dim N=2,4$
then the condition $(\dif^{D^N}\!\o_N)^{(1,2)\oplus(2,1)}=0$ is automatically
satisfied. Also, in (ii) of Proposition \ref{prop:Lic}\,, we automatically
have $\phi$ horizontally weakly conformal if $\dim N=2$\,. Therefore we have
the following result.\\

\begin{cor} \label{cor:holharmorph}
Any holomorphic map from an almost Hermitian manifold, endowed with its
Weyl connection, to a two-dimensional oriented conformal manifold is a
harmonic morphism.
\end{cor}

\indent
The Weyl connections of Example \ref{exm:minimalWeyl} and
Definition \ref{defn:canonicalWeyl} can be generalized as follows.

\begin{rem} \label{rem:PWeyl}
Let $(M,c)$ be a conformal manifold, $\dim M=m$\,, endowed with a section $P$
of ${\rm End}(TM)$ and let $\mathcal{P}=(m-1)P+P^{\,*}-\trace(P)I\!d_{TM}$\,.
It is easy to prove that if $\mathcal{P}$ is invertible at each point then for any
1-form $\a$ on $M$ there exists a unique Weyl connection $D$ on $(M,c)$ such that
$\trace_c(DP)^{\flat}=\a$\,.\\
\indent
Sufficient conditions under which $\mathcal{P}$ is invertible at each point are:\\
\indent
\quad(a) $P$ is self-adjoint and $\frac1m\trace(P)$ is not an eigenvalue of $P$
(for example, if $P=\V$ of Example \ref{exm:minimalWeyl}\,).\\
\indent
\quad(b) $P$ is skew-adjoint, invertible at each point and $m\geq3$ (for example,
if $P=J$ of Definition \ref{defn:canonicalWeyl}\,).\\
\indent
Furthermore, Proposition \ref{prop:Lic}(ii) can be easily generalized to the case
when $P$ satisfies condition (b)\,.
\end{rem}

\indent
The result of Proposition \ref{prop:Lic}(ii) can be extended as follows.

\begin{prop} \label{prop:canonicalWeyl2}
Let $(M,c_M,J^M)$ be a complex analytic almost Hermitian manifold. If $\dim M\geq4$
let $D^M$ be the Weyl connection of $(M,c_M,J^M)$\,;  if $\dim M=2$ let $D^M$ be
any Weyl connection on $(M,c_M)$\,.\\
\indent
Let $\phi:(M,c_M)\to N$ be a horizontally conformal submersion with degenerate fibres
such that ${\rm ker}\dif\!\phi$ contains $\F$ or $\widetilde{\F}$\,, where $\F$,
$\widetilde{\F}$ are the eigendistributions of $J^M$.\\
\indent
Then $\phi:(M,c_M,D^M)\to(N,D^N)$ is a harmonic map with respect to any connection $D^N$,
on $N$, and $\phi:(M,c_M,D^M)\to(N,c_N,D^N)$ is a harmonic morphism with respect to any
structure of Weyl space on $N$.\\
\indent
Conversely, if $\dim M\geq4$\,, $J^M$ is integrable and $D$ is a Weyl connection on
$(M,c_M)$ such that the foliations $\F$ and\/ $\widetilde{\F}$ are locally defined by
harmonic maps, with respect to $c_M$, $D$\,, then $D$ is the Weyl connection of
$(M,c_M,J^M)$\,.
\end{prop}
\begin{proof}
Suppose that $\F\subseteq{\rm ker}\dif\!\phi$\,. Then for any function $f$, locally
defined on $N$, the function $f\circ\phi$ is a holomorphic function of $(M,J^M)$\,.
By Remark \ref{rem:canonicalWeyl}(2)\,, $f\circ\phi$ is a harmonic function of
$(M,c_M,D^M)$ and hence $\phi:(M,c_M,D^M)\to(N,c_N,D^N)$ is a harmonic morphism with
respect to any structure of Weyl space on $N$.\\
\indent
The second statement follows from the implication (ii)$\Rightarrow$(i) of
Remark \ref{rem:canonicalWeyl}(2)\,.
\end{proof}

\indent
Let $(N^2,c_N)$ be a two-dimensional orientable conformal manifold. Then there exists
a complex structure $J^N$, uniquely determined up to sign, with respect to which
$(N^2,c_N,J^N)$ is a Hermitian manifold.\\
\indent
Let $(M^4,c_M)$ be a four-dimensional complex analytic oriented conformal manifold.
An \emph{(anti-)self-dual space at $x\in M$} is a two-dimensional vector space
$p\subseteq T_xM$ such that for some (and hence any) basis $\{X,Y\}$ of $p$ the
$2$-form $X\wedge Y$ is (anti-)self-dual; if $(M^4,c_M)$ is an oriented smooth or
real analytic manifold then an \emph{(anti-)self-dual space at $x\in M$}
is an (anti-)self-dual subspace of $(T_x^{\,\C}M,(c_M)_x^{\C})$ (see \cite{MasWoo}\,).\\
\indent
Let $(M^4,c_M)$ be a four-dimensional oriented conformal manifold endowed
with a two-dimensional nondegenerate distribution $\V$. Suppose, for simplicity, that
$\V$ (and hence, also its orthogonal complement $\H$) endowed with the conformal
structure induced by $c_M$ is orientable. We choose orientations on $(\V,c_M|_{\V})$
and $(\H,c_M|_{\H})$
so that a conformal frame $(X_1,\ldots,X_4)$ on $(M^4,c_M)$\,, adapted to the
orthogonal decomposition $TM=\V\oplus\H$, to be positive if $(X_1,X_2)$ and $(X_3,X_4)$
are positively oriented frames on $(\V,c_M|_{\V})$ and $(\H,c_M|_{\H})$\,, respectively.
Then there exists an almost complex structure $J^M$, uniquely determined up to sign,
with respect to which $(M^4,c_M,J^M)$ is a positive almost Hermitian manifold such
that $J^M(\V)=\V$. (We say that $(M^4,c_M,J^M)$ is \emph{positive} if some (and hence,
any) conformal frame of the form $(X_1,J^MX_1,X_2,J^MX_2)$ is positive;
equivalently, at some (and hence, any) point, the eigenspaces of $J^M$
are self-dual. Note that, in the smooth and real analytic categories, this just means
that $J^M$ is a positive almost complex structure on $M^4$.) It follows that the Weyl
connection of $(M^4,c_M,J^M)$ is equal to
$$D-\frac12\bigl(J^M(*_{\V}I^{\V})\bigr)^{\flat}-\frac12\bigl(J^M(*_{\H}I^{\H})\bigr)^{\flat}$$
where $D$ is the Weyl connection of $(M^4,c_M,\V)$\,, $I^{\V}$, $I^{\H}$ are the
integrability tensors of $\V$, $\H$, respectively, $*_{\V}$\,, $*_{\H}$ are the
Hodge star-operators of $(\V,c_M|_{\V})$ and $(\H,c_M|_{\H})$\,, respectively, and
$^{\flat}:TM\otimes L_M^*\!^2\longrightarrow T^*M$ is the `musical' isomorphism of
$(M^4,c_M)$ \cite{Cal-sds}\,; equivalently, the Lee form of $J^M$
with respect to any local representative $g$ of $c_M$ is equal to
$$-\trace_g(\Bv)^{\flat}-\trace_g(\Bh)^{\flat}+\bigl(J^M(*_{\V}I^{\V})\bigr)^{\flat}
+\bigl(J^M(*_{\H}I^{\H})\bigr)^{\flat}\;.$$
\indent
Let $\phi:(M^4,c_M)\to(N^2,c_N)$ be a horizontally conformal submersion, with nowhere
degenerate fibres, between oriented conformal manifolds of dimensions four and two.
Then there exists a unique
almost complex structure $J^M$ on $M^4$ with respect to which the map
$\phi:(M^4,J^M)\to(N^2,J^N)$ is holomorphic and $(M^4,c_M,J^M)$ is a positive almost
Hermitian manifold. Let $D^M$ be a Weyl connection on $(M^4,c_M)$\,.

\begin{prop}[cf.\ \cite{Woo-4d}\,] \label{prop:Woo-4d}
The following assertions are equivalent:\\
\indent
{\rm (i)} The map $\phi:(M^4,c_M,D^M)\to(N^2,c_N)$ is a harmonic morphism and
$J^M$ is integrable.\\
\indent
{\rm (ii)} The almost complex structure $J^M$ is parallel along the fibres of $\phi$\,,
with respect to $D^M$; that is, $D^M_{\;U}J^M=0$\,, $(U\in{\rm ker}\dif\!\phi)$\,.
\end{prop}
\begin{proof}
Firstly, we shall write the proof in the complex analytic category. Let $\F$ and
$\widetilde{\F}$ be the eigendistributions of $J^M$. Then $J^M$ is integrable if and
only if $\F$ and $\widetilde{\F}$ are integrable. Also, note that, assertion (ii)
holds if and only if $\F$ and $\widetilde{\F}$ are parallel, with respect to
$D^M$, along $\V\,(={\rm \ker}\dif\!\phi)$\,.\\
\indent
Let $f$ and $\widetilde{f}$ be the components of $\phi$ with respect to null local
coordinates on $(N^2,c_N)$\,. As $\phi:(M^4,c_M)\to(N^2,c_N)$ is
horizontally conformal, we have $\phi:(M^4,c_M,D^M)\to(N^2,c_N)$ is a harmonic
morphism if and only if $f$ and $\widetilde{f}$ are harmonic functions of $(M^4,c_M,D^M)$\,.
Also, we may suppose $\F\subseteq{\rm \ker}\dif\!f$\,,
$\widetilde{\F}\subseteq{\rm \ker}\dif\!\widetilde{f}$\,.\\
\indent
There exists a local frame $\bigl\{U,\widetilde{U},Y,\widetilde{Y}\bigr\}$ on $M^4$
such that $g=2(U\odot\widetilde{U}+Y\odot\widetilde{Y})$
is a local representative of
$c_M$, $U,\widetilde{U}$ are vertical, $Y,\widetilde{Y}$ are horizontal, and
$\bigl\{U,Y\bigr\}$\,, $\bigl\{\widetilde{U},\widetilde{Y}\bigr\}$ are local frames
of $\F$, $\widetilde{\F}$, respectively.\\
\indent
As $\bigl\{U,\widetilde{U},Y\bigr\}$ is a local frame of ${\rm \ker}\dif\!f$
we have $g([U,Y],Y)=0$\,. Hence, $\F$ is integrable if and only if $g([U,Y],U)=0$\,.
As $g([U,Y],U)=g(D^M_{\,U}Y,U)$ we obtain that $\F$ is integrable
if and only if $\F$ is parallel along $U$, with respect to $D^M$.\\
\indent
Also, $\trace_g(D^M\!\dif\!f)=-2g(D^M_{\,U}\widetilde{U},Y)\,\widetilde{Y}\!(f)\,.$
Hence, $f$ is a harmonic function of $(M^4,c_M,D^M)$ if and only if $\F$ is
parallel along $\widetilde{U}$, with respect to $D^M$.\\
\indent
Therefore, $\F$ is integrable and $f$ is harmonic
if and only if $\F$ is parallel, with respect to $D^M$, along $\V$.
Similarly, $\widetilde{\F}$ is integrable and $\widetilde{f}$ is harmonic if and only if
$\widetilde{\F}$ is parallel, with respect to $D^M$, along $\V$. Thus the proof is complete,
in the complex analytic category.\\
\indent
In the smooth and real analytic categories essentially the same argument applies
to the complexification $(\dif\!\phi)^{\C}:T^{\,\C}\!M\to T^{\,\C}\!N$.
\end{proof}

\section{Harmonic morphisms and twistorial maps}  \label{section:hartwist}

\indent
In this section we shall work in the complex analytic category.\\
\indent
We continue the study, initiated in the previous section, of the relation
between harmonic morphisms and twistorial maps. We start with a brief presentation
of the examples of twistorial maps with which we shall work; more details can be
found in \cite{PanWoo-sd}\,.

\begin{exm} \label{exm:twistmaps3to2}
Let $(M^3,c_M,D^M)$ be a three-dimensional Weyl space and let $(N^2,c_N)$ be a
two-dimensional conformal manifold.\\
\indent
A twistorial map $\phi:(M^3,c_M,D^M)\to(N^2,c_N)$ with nowhere degenerate fibres
is a horizontally conformal submersion whose fibres are geodesics with respect to $D^M$.
The existence of such twistorial maps is related to $(M^3,c_M,D^M)$ being
Einstein--Weyl \cite{Hit-complexmfds} (see \cite{PanWoo-sd}\,; see also
Remark \ref{rem:43to2}(1)\,, below).\\
\indent
Let $\phi:M^3\to N^2$ be a submersion with nowhere degenerate fibres
and let $p$, $\widetilde{p}$\/ be the two-dimensional degenerate distributions locally defined on $(M^3,c_M)$
such that ${\rm ker}\dif\!\phi=p\cap\widetilde{p}$\,.
Then $\phi:(M^3,c_M)\to(N^2,c_N)$ is horizontally conformal
if and only if $p$ and $\widetilde{p}$ are integrable.
It follows that $\phi:(M^3,c_M,D^M)\to(N^2,c_N)$ is twistorial if and only if
$p$ and $\widetilde{p}$
are integrable and their integral manifolds are totally-geodesic with respect
to $D^M$; note that
$\phi$ maps any such surface to a null geodesic on $(N^2,c_N)$\,.\\
\indent
By Theorems \ref{thm:Fug-Ish} and \ref{thm:BaiEel}\,,
\emph{$\phi:(M^3,c_M,D^M)\to(N^2,c_N)$ is a twistorial map if and only if it is a
harmonic morphism}.
\end{exm}

\begin{exm} \label{exm:twistmaps4to2}
Let $(M^4,c_M)$ and $(N^2,c_N)$ be oriented conformal manifolds of
dimension $4$ and $2$\,, respectively.\\
\indent
A twistorial map $\phi:(M^4,c_M)\to(N^2,c_N)$ with nowhere degenerate fibres
is a horizontally conformal submersion for which the almost complex structure
$J^M$ on $M^4$, with
respect to which $\phi:(M^4,J^M)\to(N^2,J^N)$ is holomorphic and $(M^4,c_M,J^M)$ is
a positive almost Hermitian manifold, is integrable (cf.\ \cite{Woo-4d}\,).\\
\indent
If $\phi:(M^4,c_M,D^M)\to(N^2,c_N)$ is twistorial
and $\F$\,, $\widetilde{\F}$ are the, necessarily integrable, eigendistributions of
$J^M$ then $\phi$ maps the leaves of $\F$ and $\widetilde{\F}$ to null geodesics on
$(N^2,c_N)$\,.\\
\indent
By Remark \ref{rem:canonicalWeyl}(1)\,, $\phi:(M^4,c_M)\to(N^2,c_N)$ is a
twistorial map if and only if $D^M\!J^M=0$ where $D^M$ is the Weyl connection of
$(M^4,c_M,J^M)$\,. Furthermore, if
$\phi:(M^4,c_M)\to(N^2,c_N)$ is twistorial then, by Proposition \ref{prop:Lic}\,,
$\phi:(M^4,c_M,D^M)\to(N^2,c_N)$ is a harmonic morphism. More generally,
by Proposition \ref{prop:Woo-4d}\,, \emph{if $D$ is
a Weyl connection on $(M^4,c_M)$ then
$\phi:(M^4,c_M,D)\to(N^2,c_N)$ is a twistorial harmonic morphism if and only if
$J^M$ is parallel along the fibres of $\phi$\,, with respect to $D$\,}.\\
\indent
A two-dimensional foliation $\V$ with nowhere degenerate leaves on $(M^4,c_M)$ is
twistorial if it can be locally defined by twistorial maps; note that $\V$ is
twistorial with respect to both orientations of $(M^4,c_M)$ if and only if its leaves
are totally umbilical. If $(M^4,c_M)$ is nonorientable, then $\V$ is twistorial if
its lift to the oriented $\mathbb{Z}_2$-covering space of $(M^4,c_M)$ is twistorial;
equivalently, $\V$ has totally umbilical leaves.
\end{exm}

\begin{exm} \label{exm:twistmaps4to3}
Let $(M^4,c_M)$ be an oriented four-dimensional conformal manifold and let $(N^3,c_N,D^N)$
be a three-dimensional Weyl space.\\
\indent
Let $\phi:(M^4,c_M)\to(N^3,c_N)$ be a horizontally conformal
submersion with nowhere degenerate fibres.
Let $\V={\rm ker}\dif\!\phi$\,, $\H=\V^{\perp}$ and let $D$ be the Weyl connection of
$(M^4,c_M,\V)$ (see Example \ref{exm:minimalWeyl}\,). Let ${I}^{\H}$ be the integrability
tensor of $\H$, defined by ${I}^{\H}(X,Y)=-\V[X,Y]$ for horizontal vector fields
$X$ and $Y$.\\
\indent
As both $\V$ and $\H$ are distributions of odd dimensions, the orientation of $(M^4,c_M)$
corresponds to an isomorphism between the line bundles canonically associated to
the conformal structures induced by $c_M$ on $\V$ and $\H$.
Hence, as $\V$ is one-dimensional,
both these line bundles are isomorphic to $\V$. Therefore if we apply
the Hodge star-operator $*_{\H}$ of $(\H,c_M|_{\H})$ to the integrability tensor
${I}^{\H}\in\G(\V\otimes\Lambda^2\H^*)$ of $\H$ we obtain a
horizontal one-form on $M^4$. Let $D_{\pm}$ be the Weyl partial connections on
$(\H,c_M|_{\H})$\,, over $\H$, given by $D_{\pm}=\H D\pm*_{\H}I^{\H}$
(\,\cite{Cal-sds}\,, see \cite{PanWoo-sd}\,).\\
\indent
The map $\phi:(M^4,c_M)\to(N^3,c_N,D^N)$ is twistorial, with respect
to the given orientation on $(M^4,c_M)$\,, if and only if it is
horizontally conformal and the Weyl partial connection on
$(\H,c_M|_{\H})$\,, over $\H$, determined by $D^N$ is equal to $D_{+}$\,.\\
\indent
The following assertions are equivalent for a submersion $\phi:M^4\to N^3$ with
connected nowhere degenerate fibres (\cite{Cal-sds}\,, see \cite{PanWoo-sd}\,):\\
\indent
\quad(i) There exists a Weyl connection $D^N$ on $(N^3,c_N)$ with respect to which
$\phi:(M^4,c_M)\to(N^3,c_N,D^N)$ is twistorial.\\
\indent
\quad(ii) $\phi:(M^4,c_M)\to(N^3,c_N)$ is horizontally conformal and the curvature form
of the connection induced by $D$ on $L_M$ is anti-self-dual
(that is, $\phi:(M^4,c_M)\to(N^3,c_N,D^N)$ is anti-self-dual in the sense
of \cite{Cal-sds}\,).\\
\indent
If $(M^4,c_M)$ is anti-self-dual then the following assertions can be added to this
list \cite{Cal-sds} (cf.\ \cite{Hit-complexmfds}\,, \cite{JonTod}\,; see \cite{PanWoo-sd}\,):\\
\indent
\quad(iii) There exists an Einstein--Weyl connection $D^N$ on $(N^3,c_N)$ such that
for any twistorial map $\psi$ locally defined on $(N^3,c_N,D^N)$ with values in a
conformal manifold $(P^2,c_P)$ the map $\psi\circ\phi$ from $(M^4,c_M)$ to $(P^2,c_P)$
is twistorial.\\
\indent
\quad(iv) There exists an Einstein--Weyl connection $D^N$ on $(N^3,c_N)$ such that
$\phi$ maps self-dual surfaces on $(M^4,c_M)$ to degenerate surfaces
on $(N^3,c_N)$ which are totally-geodesic with respect to $D^N$.\\
\indent
It follows that if $\phi:(M^4,c_M)\to(N^3,c_N,D^N)$ is twistorial then
$(M^4,c_M)$ is anti-self-dual if and only if $(N^3,c_N,D^N)$ is Einstein--Weyl
\cite{Cal-sds} (cf.\ \cite{Hit-complexmfds}\,, \cite{JonTod}\,; see \cite{PanWoo-sd}\,).
Furthermore, if $(M^4,c_M)$ is anti-self-dual then, locally, $\phi$ corresponds
to a submersion $Z(\phi)$ from the twistor space $Z(M)$ of $(M^4,c_M)$ onto the
twistor space $Z(N)$ of $(N^3,c_N,D^N)$ which maps each twistor line on $Z(M)$
diffeomorphically onto a twistor line on $Z(N)$\,; the map $Z(\phi)$ is the
\emph{twistorial representation} of $\phi$\,.\\
\indent
As in Example \ref{exm:twistmaps4to2}\,, a one-dimensional foliation $\V$ with nowhere
degenerate leaves is twistorial if it can be locally defined by twistorial maps.
Note that $\V$ is twistorial with respect to both orientations of $(M^4,c_M)$ if and
only if it is locally generated by conformal vector fields \cite{Cal-sds}
(this follows from the equivalence (i)$\iff$(ii)\,, above, and
Example \ref{exm:minimalWeyl}\,). If $(M^4,c_M)$ is nonorientable then $\V$ is twistorial
if its lift to the oriented $\mathbb{Z}_2$-covering space of $(M^4,c_M)$ is twistorial;
equivalently, $\V$ is locally generated by nowhere zero conformal vector fields.
\end{exm}

\indent
Next, we give necessary and sufficient conditions under which a map
between Weyl spaces of dimensions four and three is a twistorial
harmonic morphism.

\begin{thm} \label{thm:twistharmorph4to3}
Let $(M^4,c_M,D^M)$ be an oriented Weyl space of dimension four and let
$(N^3,c_N)$ be a conformal manifold of dimension three. We denote by
$\tD^M\subseteq TP_M$ the connection induced by $D^M$ on
the bundle $\p_M:P_M\to M$ of self-dual spaces on $(M^4,c_M)$\,.
Also, we denote by $\p_N:P_N\to N$ the bundle of two-dimensional degenerate
spaces on $(N^3,c_N)$\,.\\
\indent
Let $\phi:(M^4,c_M)\to(N^3,c_N)$ be a horizontally conformal submersion
with connected nowhere degenerate fibres and let $D$ be the Weyl connection
of $(M^4,c_M,\V)$\,. We denote by $\varPhi:P_M\to P_N$ the bundle
map defined by $\varPhi(p)=\dif\!\phi(p)$\,, $(p\in P_M)$\,.\\
\indent
{\rm (a)} Let $D_N$ be a Weyl connection on $(N^3,c_N)$\,. Any two of the
following assertions imply the third:\\
\indent
\quad{\rm (a1)} $\phi:(M^4,c_M,D^M)\to(N^3,c_N,D^N)$ is a harmonic morphism.\\
\indent
\quad{\rm (a2)} $\phi:(M^4,c_M)\to(N^3,c_N,D^N)$ is twistorial.\\
\indent
\quad{\rm (a3)} The fibres of $\varPhi$ are tangent to $\tD^M$.\\
\indent
{\rm (b)} The following assertions are equivalent:\\
\indent
\quad{\rm (b1)} There exists a Weyl connection $D^N$ on $(N^3,c_N)$ and a section
$k$ of the dual of the line bundle $L_N$ of $(N^3,c_N)$ such that
$\phi:(M^4,c_M,D^M)\to(N^3,c_N,D^N)$ is a twistorial harmonic morphism,
and the vertical component of $D^M-D$ is equal to $\tfrac12\,k$\,, under the
isomorphism $\V=\phi^*(L_N)$ corresponding to the orientation of $(M^4,c_M)$.\\
\indent
\quad{\rm (b2)} $\tD^M$ is projectable, with respect to $\varPhi$\,, onto a three-dimensional
distribution on $P_N$.\\
\indent
\quad{\rm (b3)} There exists a Weyl connection $D^N$ on $(N^3,c_N)$ and a section
$k$ of $L_N^*$ such that $\dif\!\varPhi(\tD^M)=\tD^{\nabla}$ where
$\tD^{\nabla}\subseteq TP_N$ is the connection
on $P_N$ induced by the connection $\nabla$ on $L_N^*\otimes TN$ defined by
\begin{equation} \label{e:nabla}
\nabla_X\xi=D^N_{\,X}\xi+\tfrac12\,kX\times\xi
\end{equation}
for any local sections $X$ of $TN$ and $\xi$ of $L_N^*\otimes TN$.\\
\indent
\quad{}Furthermore, if assertion {\rm (b1)}\,,\,{\rm (b2)}\,,\,{\rm (b3)} holds then
$D^N$ and $k$ are determined by
\begin{equation} \label{e:ktwistharmorph4to3}
D^M=D+\tfrac12\,\bigl(k+*_{\H}\!I^{\H}\,\bigr)\quad \textit{and}\quad D^N=D_+\;.
\end{equation}
\end{thm}
\begin{proof}
To prove (a)\,, we claim that (a3) is equivalent to the following equality of
partial connections, over $\H$,
\begin{equation} \label{e:a3}
\H D^M=\H D+\tfrac12*_{\H}\!I^{\H}\;.
\end{equation}
Indeed, note that (a3) holds if and only if for any positively
oriented conformal local frame $\bigl(X_0,X_1,X_2,X_3\bigr)$ such that $X_0$ is
tangent to the fibres of $\phi$ and $X_1$\,,\,$X_2$\,,\,$X_3$ are basic we have
\begin{equation} \label{e:a3frame}
g\bigl(D^M_{\,X_0}(X_0+{\rm i}X_1),X_2+{\rm i}X_3\bigr)=0
\end{equation}
where $g$ is the local representative of $c_M$ determined by
$\bigl(X_0,X_1,X_2,X_3\bigr)$\,. An easy calculation shows that \eqref{e:a3frame} is
equivalent to
$$\bigl(\a_M-\bigl(\trace_g(\Bv)^{\flat}+\tfrac12*_{\H}\!I^{\H}\bigr)\bigr)(X_2+{\rm i}X_3)=0$$
where $\a_M$ is the Lee form of $D^M$ with respect to $g$\,. Now, the equivalence
(a3)$\iff$\eqref{e:a3} follows easily.\\
\indent
Assertion (a) follows from Proposition \ref{prop:fundamn} and
Example \ref{exm:twistmaps4to3}\,.\\
\indent
To prove (b)\,, firstly note that (b3)$\Longrightarrow$(b2) is trivial. Conversely,
if (b2) holds then, as $\p_N$ is proper, $\dif\!\varPhi(\tD^M)$ is a connection on
$P_N$ whose holonomy group is contained in the group of (complex analytic) diffeomorphisms
of $\C\!P^1$\,. Hence, the holonomy group of $\dif\!\varPhi(\tD^M)$ is contained
in $PGL(2,\C)\bigl(=SL(2,\C)/\mathbb{Z}_2\bigr)$\,. As $SO(3,\C)$ is the adjoint
group of $SL(2,\C)$\,, we have $SO(3,\C)=SL(2,\C)/\mathbb{Z}_2$\,. Therefore any
connection on $P_N$
corresponds to a connection on the oriented Riemannian bundle $(L_N^*\otimes TN,c_N)$\,.
In particular, there exists a unique connection $\nabla$ on $(L_N^*\otimes TN,c_N)$
which induces the connection $\dif\!\varPhi(\tD^M)$ on $P_N$.
Note that $\nabla$ canonically determines a projective structure on $N^3$ which
has the property that any of its geodesics which is null at some point is null everywhere.
Furthermore, a null curve on $(N^3,c_N)$ is a geodesic of $\nabla$ if and only if
its velocity vector field $Y$ has the property that $Y^{\perp}$ is parallel with
respect to $\nabla$. {}From the fact that $\phi$ is horizontally conformal, it follows
that $\nabla$ and $c_N$ have the same null geodesics;
equivalently, there exists a Weyl connection $D^N$ on $(N^3,c_N)$ and a section
$k$ of $L_N^*$ so that \eqref{e:nabla} holds. Thus, we have proved that
(b2)$\Longrightarrow$(b3)\,.\\
\indent
It is obvious that if (b3) holds then $\nabla$ and $D^N$ determine the same projective
structure. It follows that $\phi:(M^4,c_M)\to(N^3,c_N,D^N)$ is twistorial.
Thus, by (a)\,, if (b3) holds then
$\phi:(M^4,c_M,D^M)\to(N^3,c_N,D^N)$ is a twistorial harmonic morphism.
Further, as $\phi:(M^4,c_M)\to(N^3,c_N,D^N)$ is twistorial,
the partial connection over $\H$ determined by the pull-back by $\phi$
of $D^N$ is equal to $D_+$\,. Also, an argument similar to the proof of (a)
shows that a field $p$ of self-dual spaces, over a horizontal curve on $(M^4,c_M)$\,,
is parallel with respect to $D^M$ if and only if $\H p$ is parallel with respect to the
partial connection $\widehat{\nabla}$ on $\V^*\!\otimes\!\H$, over $\H$, defined by
$$\widehat{\nabla}_X\xi=(D_+)_X\xi+\tfrac12\,\widehat{k}X\times\xi$$ for any local
sections $X$ of $\H$ and $\xi$ of $\V^*\!\otimes\!\H$, where $\widehat{k}$
is the section of $\V^*$ defined by
$$D^M=D+\tfrac12\,\bigl(\,\widehat{k}+*_{\H}I^{\H}\,\bigr)\;.$$
Thus, the partial connection over $\H$ induced by the pull-back by $\phi$
of $\nabla$ is equal to $\widehat{\nabla}$. Consequently, $\widehat{k}$ is the
pull-back by $\phi$ of $k$\,. In particular, (b3)$\Longrightarrow$(b1)\,.\\
\indent
We have proved that (b3) is equivalent to the existence of a Weyl connection
$D^N$ on $(N^3,c_N)$ and a basic section $k$
of $\V^*\bigl(=\phi^*(L_N^*)\bigr)$ such that \eqref{e:ktwistharmorph4to3} holds.
On the other hand, by the proof of (a)\,, the map
$\phi:(M^4,c_M,D^M)\to(N^3,c_N,D^N)$ is a twistorial harmonic morphism, for some
Weyl connection $D^N$ on $(N^3,c_N)$\,, if and only if there exists a
section $k$ of $\V^*$ such that \eqref{e:ktwistharmorph4to3} holds.
Thus, (b1)$\Longrightarrow$(b3).\\
\indent
The theorem is proved.
\end{proof}

\begin{rem} \label{rem:4to3geodsfibres}
Let $\phi:(M^4,c_M,D^M)\to(N^3,c_N,D^N)$ be a twistorial
harmonic morphism with nowhere degenerate fibres.
{}From Theorems \ref{thm:BaiEel} and \ref{thm:twistharmorph4to3}(a)
it follows that the fibres of $\phi$ are
geodesics with respect to $D^M$ if and only if $\H$ is integrable;
in particular, $\phi$ is also twistorial with respect to the
reversed orientation of $(M^4,c_M)$\,, and hence, the fibres of
$\phi$ are locally generated by nowhere zero conformal vector
fields.
\end{rem}

\indent
For later use, we formulate the following definition (cf.\ \cite{GauTod}\,,
\cite{CalPed}\,).

\begin{defn} \label{defn:GauTod}
Let $(N^3,c,D)$ be a three-dimensional Weyl space. We say that $(N^3,c,D)$ is
a \emph{Gauduchon-Tod space} if there exist three one-dimensional conformal foliations
by geodesics $\V_1$\,,\,$\V_2$\,,\,$\V_3$ which are orthogonal on each other
and such that $\V_j^{\perp}$, considered with the conformal structure
induced by $c$\,, is orientable, $j=1,2,3$\,.
\end{defn}

\begin{rem} \label{rem:GauTod}
Let $(N^3,c)$ be a three-dimensional conformal manifold and let $L$ be the
associated line bundle. Then, it is known that, locally,
the following assertions are equivalent (the equivalences (i)$\iff$(ii)$\iff$(iii)
are proved in \cite{GauTod} whilst (i)$\iff$(iv)$\iff$(v) follow from
\cite{Hit-complexmfds}\,):\\
\indent
(i) There exists a Weyl connection $D$ on $(N^3,c)$ with respect to which
$(N^3,c,D)$ is Gauduchon-Tod.\\
\indent
(ii) There exists an Einstein--Weyl connection $D$ on $(N^3,c)$ whose scalar
curvature $s^D$ is given by $s^D=\tfrac32\,k^2$ for a section $k$ of $L^*$
which satisfies $*_{c\,}Dk=F^D$ where $F^D$ is the curvature form of the connection,
on $L$\,, corresponding to $D$\,.\\
\indent
(iii) There exists a Weyl connection $D$ on $(N^3,c)$ and a section
$k$ of $L^*$ such that the connection $\nabla$ on $L^*\otimes TN$, defined by
\begin{equation*}
\nabla_X\xi=D_X\xi+\tfrac12\,kX\times\xi
\end{equation*}
for any local sections $X$ of $TN$ and $\xi$ of $L^*\otimes TN$, is flat.\\
\indent
(iv) There exists a flat connection on $(L^*\otimes TN,c)$\,.\\
\indent
(v) There exists an Einstein--Weyl connection $D$ on $(N^3,c)$ for which there exists
a submersion from its twistor space onto $\C\!P^1$ whose fibres are transversal to the
twistor lines.
\end{rem}

\indent
We end this section with the following result.

\begin{cor} \label{cor:ObataGauTodharmorph}
Let $\phi:(M^4,c_M,D^M)\to(N^3,c_N,D^N)$ be a twistorial harmonic morphism with
nowhere degenerate fibres from an oriented Weyl space of dimension four to a
Weyl space of dimension three; let $k\in\G(\V^*)$ be the vertical component
of $2(D^M-D)$\,. Then, locally, the following assertions are equivalent:\\
\indent
\quad{\rm (i)} $D^M$ is the Obata connection of a hyper-Hermitian
structure on $(M^4,c_M)$\,.\\
\indent
\quad{\rm (ii)} $(N^3,c_N,D^N)$ is Gauduchon-Tod and $*_ND^Nk=F^{D^N}$.\\
\indent
\quad{\rm (iii)} $(M^4,c_M,D^M)$ is Einstein-Weyl anti-self-dual and $k$ is basic.\\
\indent
Conversely, let $\phi:(M^4,c_M)\to(N^3,c_N,D^N)$ be a twistorial map,
with nowhere degenerate fibres, from an oriented conformal manifold of dimension
four to a Gauduchon-Tod space, and let $D^M$ be the Obata connection of the hyper-Hermitian
structure induced on $(M^4,c_M)$\,. Then $\phi:(M^4,c_M,D^M)\to(N^3,c_N,D^N)$
is a harmonic morphism.
\end{cor}
\begin{proof}
With the same notations as in Theorem \ref{thm:twistharmorph4to3}\,, let
$\F\subseteq\tD^M$ be the twistor distribution defined by
$\F_p\subseteq T_pP_M$ is the horizontal lift of
$p\subseteq T_{\p_M(p)}M$, with respect to $D^M$, ($p\in P_M$)\,.
Then, $(M^4,c_M)$ is anti-self-dual if and only if $\F$ is integrable.
Also, $(M^4,c_M,D^M)$ is Einstein--Weyl anti-self-dual if and only if
$\tD^M$ is projectable with respect to $\F$ (that is, $[X,Y]$ is a local section
of $\tD^M$ if $X$ and $Y$ are local sections of $\tD^M$ and $\F$, respectively).\\
\indent
As (i) is satisfied if and only if the connection $\tD^M$ is trivial,
the equivalence (i)$\iff$(ii) follows easily from Theorem \ref{thm:twistharmorph4to3}(a)
and Remark \ref{rem:GauTod}\,.\\
\indent
If (i) holds, from Theorem \ref{thm:twistharmorph4to3}(a) we obtain that
$\tD^M$ is projectable with respect $\varPhi$. Thus, by
Theorem \ref{thm:twistharmorph4to3}(b)\,, $k$ is basic.
Hence, (i)$\Longrightarrow$(iii)\,.\\
\indent
If (iii) holds then, by Theorem \ref{thm:twistharmorph4to3}(b)\,,
$\tD^M$ is projectable with respect to $\varPhi$. It follows easily that $\tD^M$
is integrable. Thus, locally, (iii)$\Longrightarrow$(i)\,.\\
\indent
The converse statement follows from Theorem \ref{thm:twistharmorph4to3}(b)\,.
\end{proof}

\section{Relations between the twistoriality of harmonic morphisms and the Ricci tensor} \label{section:Riccitwistsharmorphs}

\indent
In this section we shall work in the complex analytic category. Here, by a
\emph{real conformal manifold} we mean the (germ-unique) complexification of a real
analytic conformal manifold; similarly, we sometimes work with \emph{real Weyl spaces}.\\
\indent
The following lemma follows from a straightforward computation.

\begin{lem}[cf.\ \cite{BaiWoo1}\,] \label{lem:RicciY}
Let $(M,c,D)$ be a Weyl space, $\dim M=3,4$\,, and let ${\rm Ric}$ be its Ricci tensor.
Let $\F$ be a foliation by null geodesics on $(M,c,D)$ such that $\F^{\perp}$ is integrable.\\
\indent
{\rm (i)} If $\dim M=3$ then
\begin{equation*}
\Ric(Y,Y)=Y(g(D_UU,Y))-g(D_UU,Y)^2
\end{equation*}
where $\bigl\{U,Y,\widetilde{Y}\bigr\}$ is a local frame on $M$ such that $Y$ is a local
section of $\F$\,, $D_YY=0$ and
$g=U\odot U+2\,Y\odot\widetilde{Y}$ is a local representative of $c$\,.\\
\indent
{\rm (ii)} If $\dim M=4$ then
\begin{equation*}
\Ric(Y,Y)=2\bigl[Y(g(D_U\widetilde{U},Y))-g(D_U\widetilde{U},Y)^2-g([U,Y],U)g([\widetilde{U},Y],\widetilde{U})\bigr]
\end{equation*}
where $\bigl\{U,\widetilde{U},Y,\widetilde{Y}\bigr\}$ is a local frame on $M$ such that $Y$
is a local section of $\F$\,, $D_YY=0$ and
$g=2(U\odot\widetilde{U}+Y\odot\widetilde{Y})$ is a local representative of\/ $c$\,. \qed
\end{lem}

\begin{rem}
1) In Lemma \ref{lem:RicciY}(i) the condition $\F^{\perp}$ integrable is superfluous.
It follows that from any three-dimensional conformal manifold we can, locally, define
horizontally conformal submersions with one-dimensional nowhere degenerate fibres tangent
to any given direction at a point. A similar statement holds for real analytic
three-dimensional conformal manifolds.\\
\indent
2) A relation slightly longer than in Lemma \ref{lem:RicciY}(ii) can be obtained,
for a foliation $\F$ by null geodesics on a four-dimensional Weyl space, without
the assumption $\F^{\perp}$ is integrable.\\
\end{rem}

\begin{prop}[cf.\ \cite{BaiWoo1}\,,\,\cite{Woo-4d}\,]
\label{prop:43to2}
Let $(M,c_M,D^M)$ be a Weyl space, and let
$(N,c_N)$ be a conformal manifold, $\dim M=3\,,4$\,, $\dim N=2$\,.
If $\dim M=4$ suppose that $\phi$ is real (that is, $\phi$ is the
(germ-unique) complexification of a real analytic map and, in
particular, $(M,c_M,D^M)$ and $(N,c_N)$ are complexifications of a
real analytic Weyl space and of a real analytic conformal manifold,
respectively).\\
\indent
Let $\phi:(M,c_M)\to(N,c_N)$ be a horizontally conformal submersion with nowhere
degenerate fibres; suppose $\trace_{c_M}(D\!\dif\!\phi)=0$
along a hypersurface transversal to the fibres of $\phi$\,.\\
\indent
Then any two of the following assertions imply the third:\\
\indent
\quad{\rm (i)} $\phi:(M,c_M,D^M)\to(N,c_N)$ is a harmonic morphism.\\
\indent
\quad{\rm (ii)} $\phi:(M,c_M,D^M)\to(N,c_N)$ is twistorial.\\
\indent
\quad{\rm (iii)} The trace free symmetric part of the horizontal component of the Ricci
tensor of $D^M$ is zero.
\end{prop}
\begin{proof}
If $\dim M=3$ then (i)$\iff$(ii)\,. Also, we can find a local frame $\bigl\{U,Y,\widetilde{Y}\bigr\}$\,,
as in Lemma \ref{lem:RicciY}(i)\,, such that $U$ is tangent to the fibres of $\phi$\,.
Then assertion (ii) is equivalent to $g(D^M_{\;U}U,Y)=g(D^M_{\;U}U,\widetilde{Y})=0$\,.
On the other hand, assertion (iii) is equivalent to
$\RicM(Y,Y)=\RicM(\widetilde{Y},\widetilde{Y})=0$ where
$\RicM$ is the Ricci tensor of $D^M$\,. Thus, if $\dim M=3$\,, the proof follows
from Lemma \ref{lem:RicciY}(i)\,.\\
\indent
Suppose $\dim M=4$\,. Then we can find a local frame
$\bigl\{U,\widetilde{U},Y,\widetilde{Y}\bigr\}$
like in Lemma \ref{lem:RicciY}(ii)\,, such that $U$ and $\widetilde{U}$ are tangent to the
fibres of $\phi$\,. Moreover, we may assume $g$ oriented such that $\F_+=Span(U,Y)$ and
$\widetilde{\F}_+=Span(\widetilde{U},\widetilde{Y})$ are self-dual whilst
$\F_-=Span(\widetilde{U},Y)$ and $\widetilde{\F}_-=Span(U,\widetilde{Y})$ are anti-self-dual.
Then assertion (ii) is equivalent to the fact that either $\F_+$\,, $\widetilde{\F}_+$
are integrable or $\F_-$\,, $\widetilde{\F}_-$ are integrable. On the other
hand, assertion (i) is equivalent to
$g(D^M_{\,U}\widetilde{U},Y)=g(D^M_{\,U}\widetilde{U},\widetilde{Y})=0$ (see the proof
of Proposition \ref{prop:Woo-4d}\,). Thus, if $\dim M=4$\,, the proof follows
from Lemma \ref{lem:RicciY}(ii)\,.
\end{proof}

\begin{rem} \label{rem:43to2}
1) If $\dim M=3$ then assertions (i)\,, (ii)\,, (iii) of
Proposition \ref{prop:43to2} are equivalent.\\
\indent
It follows that if $(M,c_M,D^M)$ is a three-dimensional Weyl space from which can
be locally defined more than $k=6$ harmonic morphisms with one-dimensional nowhere
degenerate fibres then $(M,c_M,D^M)$ is Einstein--Weyl; in the smooth category,
the same statement holds with $k=2$ \cite{CalPed} (cf.\ \cite{BaiWoo1}\,).\\
\indent
2) In the smooth category, suppose assertion (i) of Proposition \ref{prop:43to2} holds.\\
\indent
\quad(a) If $\dim M=3$ then assertion (iii) also holds (cf.\ \cite{BaiWoo1}\,).\\
\indent
\quad(b) If $\dim M=4$ then the implication (ii)$\Rightarrow$(iii) holds on $M$
whilst the implication (iii)$\Rightarrow$(ii) holds locally on a dense open set of $M$
(cf.\ \cite{Woo-4d}\,).\\
\indent
3) Proposition \ref{prop:43to2} also holds for any horizontally conformal submersion
$\phi:(M,c_M,D^M)\to(N,c_N,D^N)$ with nowhere degenerate fibres from a Weyl space
to an Einstein--Weyl space, $\dim M=4$\,, $\dim N=3$\,
(see Proposition \ref{prop:4to23}\,, below, for an extension of this fact).\\
\indent
4) For the implications (i)\,,\,(ii)$\Rightarrow$(iii) and
(ii)\,,\,(iii)$\Rightarrow$(i)\,, of Proposition \ref{prop:43to2}\,, it is not necessary
to assume that $\phi$ is real when $\dim M=4$.
\end{rem}

\indent
The proof of the following lemma is omitted.

\begin{lem}[cf.\ \cite{PanWoo-exm}\,] \label{lem:RicciY4to3}
Let $\phi:(M,c_M,D^M)\to(N,c_N,D^N)$ be a submersive harmonic morphism with nowhere
degenerate fibres between Weyl spaces, $\dim M=4$\,, $\dim N=3$\,.\\
\indent
Let $A_{\pm}=D_{\pm}-D^N$. Then for any horizontal null vector $Y$ we have
\begin{equation*}
\RicM(Y,Y)-\RicN(\dif\!\phi(Y),\dif\!\phi(Y))=-\,\tfrac12\,A_+(Y)\,A_-(Y)
\end{equation*}
where $\RicM$ and $\RicN$ are the Ricci tensors of $D^M$ and $D^N$\,,
respectively.   \qed
\end{lem}

\indent
The following result follows from Lemmas \ref{lem:RicciY}(ii) and \ref{lem:RicciY4to3}\,.

\begin{prop}[cf.\ \cite{PanWoo-exm}\,] \label{prop:4to23}
Let $\phi:(M,c_M,D^M)\to(N,c_N,D^N)$ be a nonconstant harmonic morphism with nowhere
degenerate fibres between Weyl spaces, $\dim M=4$\,, $\dim N=2,\,3$\,.
If $\dim N=2$ suppose that $\phi$ is real.\\
\indent
Let $\RicM$ and $\RicN$ be the Ricci tensors of $D^M$ and $D^N$\,, respectively.
Then the following assertions are equivalent:\\
\indent
\;{\rm (i)} $\phi$ is twistorial.\\
\indent
\;{\rm (ii)} The trace free symmetric part of the horizontal component of
$\RicM-\phi^*\bigl(\RicN\bigr)$ is zero.   \qed
\end{prop}

\section{Harmonic morphisms from four-dimensional Einstein--Weyl spaces}     \label{section:4to23}

\indent
In this section, we shall work in the complex analytic category; here, as in the previous
section, by a real conformal manifold (Weyl space) we mean the complexification of a real
analytic conformal manifold (Weyl space).\\
\indent
{}From Proposition \ref{prop:4to23}\,, we obtain the following two corollaries.

\begin{cor}[cf.\ \cite{Woo-4d}\,,\,\cite{PanWoo-sd}\,] \label{cor:4to2}
Let $\phi:(M^4,c_M,D^M)\to(N^2,c_N)$ be a real submersive harmonic morphism
from a four-dimensional real Einstein--Weyl space to a two-dimensional real conformal manifold.\\
\indent
Then $\phi:(M^4,c_M)\to(N^2,c_N)$ is a twistorial map. Furthermore, if $(M^4,c_M)$
is nonorientable then the fibres of $\phi$ are totally geodesic. \qed
\end{cor}

\begin{cor}[cf.\ \cite{Pan-4to3}\,,\,\cite{PanWoo-exm}\,,\,\cite{PanWoo-sd}\,]  \label{cor:4to3}
Let $(M^4,c_M,D^M)$ be an Einstein--Weyl space of dimension four and
let $\phi:(M^4,c_M,D^M)\to(N^3,c_N,D^N)$ be a submersive harmonic morphism with nowhere
degenerate fibres to a Weyl space of dimension three.\\
\indent
Then $\phi:(M^4,c_M)\to(N^3,c_N,D^N)$ is a twistorial map if and only if $D^N$ is an
Einstein--Weyl connection. \qed
\end{cor}

\begin{rem}
We do not know whether or not it is true that if $(M^4,c_M,D^M)$ is Einstein--Weyl
and $\phi:(M^4,c_M,D^M)\to(N^3,c_N,D^N)$ is a submersive harmonic morphism with
nowhere degenerate fibres then, locally, there exists a Weyl connection $D$ on
$(N^3,c_N)$ with respect to which $\phi:(M^4,c_M)\to(N^3,c_N,D)$ is twistorial.\\
\indent
{}From \cite{Pan-4to3}\,, \cite{PanWoo-sd}\,, it follows that this holds if $D^M$ and $D^N$
are Levi-Civita connections of representatives of $c_M$ and $c_N$, respectively.
Furthermore, if we also assume $D\neq D^N$ then $\phi$ must be of Killing type;
for example, the Killing vector field
$V=x_1\,\partial/\partial{x_2}-x_2\,\partial/\partial{x_1}$
on $\C^{\!4}$ defines such harmonic morphisms (this follows from
Example \ref{exm:twistmaps4to3} by noting that $V^{\perp}$ is integrable but
the nondegenerate orbits of $V$ are nowhere geodesic).
\end{rem}

\indent
Next, we give the necessary and sufficient conditions under which on a real
Einstein--Weyl space of dimension four there can be defined, locally, at least
five distinct real foliations of dimension two which produce harmonic morphisms
(cf.\ Remark \ref{rem:43to2}(1)\,).

\begin{thm} \label{thm:five4to2}
Let $(M^4,c_M,D^M)$ be an orientable real Einstein--Weyl space of dimension four.\\
\indent
Then, locally, there can be defined on $(M^4,c_M,D^M)$ \emph{five} distinct real
foliations of dimension two which produce harmonic morphisms if and only if
one of the following two assertions holds:\\
\indent
\quad{\rm (i)} $D^M$ is the Weyl connection of a Hermitian structure locally defined
on $(M^4,c_M)$;\\
\indent
\quad{\rm (ii)} $(M^4,c_M)$ is anti-self-dual, with respect to a suitable orientation,
and $D^M$ is the Levi-Civita connection of a local Einstein representative of $c_M$\,.
\end{thm}
\begin{proof}
Let $\V_k$\,, $k=1\,,\ldots,\,5$\,, be five
distinct two-dimensional foliations which produce harmonic morphisms on
$(M^4,c_M,D^M)$\,. By passing to a conformal covering, if necessary, we may assume
that $\V_k$\,, $k=1\,,\ldots,\,5$\,,
endowed with the conformal structures induced by $c_M$\,, are orientable.
Also, as $M^4$ is locally compact, we may suppose that
$\V_k(x)$\,, $k=1\,,\ldots,\,5$\,, are distinct at each $x\in M$.\\
\indent
By Corollary \ref{cor:4to2}\,, for each $k=1\,,\ldots,\,5$\,, we have
$\V_k$ twistorial, with respect to a suitable orientation on $(M^4,c_M)$\,. Thus,
there exists an orientation on $(M^4,c_M)$ with respect to which, after a
renumerotation, $\V_k$ is twistorial, for $k=1\,,\,2\,,\,3$\,.\\
\indent
Let $J_k$ be the positive Hermitian structures on
$(M^4,c_M)$\,, uniquely determined up to sign, such that $J_k(\V_k)=\V_k$\,,
$k=1\,,\,2\,,\,3$\,.
Then either there exists $k\in\{1\,,\,2\,,\,3\}$ such that $D^M$ is the Weyl
connection of $(M^4,c_M,J_k)$ or for any $k\in\{1\,,\,2\,,\,3\}$ we have
$D^MJ_k\neq0$\,. In the
latter case, from Proposition \ref{prop:Woo-4d} it follows that, for
$k=1\,,\,2\,,\,3$\,, the complex structure $J_k$ determines $\V_k$
(cf.\ \cite{Woo-4d}\,) and hence we obtain that $J_k\pm J_l\neq0$ for
$k\neq l$\,; thus, $(M^4,c_M)$ is anti-self-dual (see \cite{Sal}\,).
Now, as $(M^4,c_M,D^M)$ is
anti-self-dual and Einstein--Weyl, by a result of H.~Pedersen, A.~Swann, and
D.M.J.~Calderbank (see \cite{Cal-F} and the references therein), we have that, locally,
either $D^M$ is the Obata connection of a hyper-Hermitian structure on $(M^4,c_M)$
or $D^M$ is the Levi-Civita connection of an Einstein representative of $c_M$\,.\\
\indent Conversely, if (i) or (ii) holds then, locally, there can be
defined infinitely many two-dimensional foliations on
$(M^4,c_M,D^M)$ which produce harmonic morphisms. Indeed, if
$(M^4,c_M,D^M)$ satisfies assertion (i) then, this follows easily
from Remark \ref{rem:canonicalWeyl}(ii) whilst, if $(M^4,c_M,D^M)$
satisfies assertion (ii) then, this is, essentially, proved in
\cite{Woo-4d}\,. The theorem is proved.
\end{proof}

\begin{rem}
Let $(M^4,c_M,D^M)$ be an orientable four-dimensional Einstein--Weyl space endowed
either with \emph{nine} distinct foliations, by degenerate hypersurfaces, which produce
harmonic morphisms or with \emph{thirteen} distinct foliations, by nondegenerate
surfaces, which produce harmonic morphisms.\\
\indent
Then, arguments as in the proof of Theorem \ref{thm:five4to2} show that, up to conformal
covering spaces, one of the following two assertions holds:\\
\indent
\quad(i) The connection induced by $D^M$ on the bundle of self-dual spaces on
$(M^4,c_M)$\,, considered with a suitable orientation, admits a reduction to the group
of affine transformations of $\C$.\\
\indent
\quad(ii) $(M^4,c_M)$ is anti-self-dual, with respect to a suitable orientation,
and $D^M$ is the Levi-Civita connection of an Einstein representative of $c_M$\,.
\end{rem}

\indent
Finally, we describe the harmonic morphisms with nowhere degenerate fibres between
Einstein--Weyl spaces of dimensions four and three.

\begin{thm} \label{thm:4to3}
Let $(M^4,c_M,D^M)$ and $(N^3,c_N,D^N)$ be Einstein--Weyl spaces of dimension four and
three, respectively, and let $\phi:M^4\to N^3$ be a submersion with nowhere degenerate
fibres.\\
\indent
If $(M^4,c_M)$ is orientable then $\phi:(M^4,c_M,D^M)\to(N^3,c_N,D^N)$ is a
harmonic morphism if and only if $(M^4,c_M)$ is anti-self-dual, with respect to a
suitable orientation, $\phi:(M^4,c_M)\to(N^3,c_N,D^N)$ is twistorial
and, locally, either\\
\indent
\quad{\rm (i)} $(N^3,c_N,D^N)$ is Gauduchon-Tod and $D^M$ is the Obata connection
of the hyper-Hermitian structure induced on $(M^4,c_M)$ or\\
\indent
\quad{\rm (ii)} $D^M$ is the Levi-Civita connection of an Einstein representative $g$
of $c_M$ with nonzero scalar curvature and the fibres of the twistorial
representation $Z(\phi):Z(M)\to Z(N)$ of $\phi$ are tangent to the contact distribution
on $Z(M)$ corresponding to $g$\,.\\
\indent
If $(M^4,c_M)$ is nonorientable then $\phi:(M^4,c_M,D^M)\to(N^3,c_N,D^N)$ is a harmonic
morphism if and only if, locally, $D^M$ and $D^N$ are the Levi-Civita
connections of constant curvature representatives of $c_M$ and $c_N$, respectively,
and $\phi$ is a harmonic morphism of warped product type.
\end{thm}
\begin{proof}
\indent
If $\phi:(M^4,c_M)\to(N^3,c_N,D^N)$ is twistorial and (i) or (ii) holds
then (a2) and (a3) of Theorem \ref{thm:twistharmorph4to3} are satisfied and hence
$\phi:(M^4,c_M,D^M)\to(N^3,c_N,D^N)$ is a harmonic morphism.\\
\indent
Conversely, if $\phi:(M^4,c_M,D^M)\to(N^3,c_N,D^N)$ is a harmonic morphism then,
by Corollary \ref{cor:4to3}\,, we have that $\phi:(M^4,c_M)\to(N^3,c_N,D^N)$
is a twistorial map. {}From \cite{Cal-sds} (see \cite{PanWoo-sd}\,), it follows
that either $(M^4,c_M)$ is orientable and anti-self-dual, with respect to a suitable
orientation, or $(M^4,c_M)$ is flat.\\
\indent
If $(M^4,c_M)$ is orientable, as $(M^4,c_M,D^M)$ is Einstein--Weyl, then, as in the
proof of Theorem \ref{thm:five4to2}\,, we have that, locally,
either $D^M$ is the Obata connection of a hyper-Hermitian structure on $(M^4,c_M)$
or $D^M$ is the Levi-Civita connection of an Einstein representative of $c_M$\,;
if $(M^4,c_M)$ is nonorientable then it must be flat and hence, by a result of
M.G.~Eastwood and K.P.~Tod, we have that, locally, $D^M$ is the  Levi-Civita connection
of some constant curvature representative of $c_M$ (see \cite{Cal-F} and the references
therein, and note that similar calculations prove that these two results also hold in
the complex analytic category).\\
\indent
If $(M^4,c_M)$ is orientable then, from Theorem \ref{thm:twistharmorph4to3}(a)
it follows that, locally, we have the alternative (i) or (ii)\,.\\
\indent
If $(M^4,c_M)$ is nonorientable then, by passing to its oriented $\mathbb{Z}_2$-covering,
we obtain that $\phi$ has integrable horizontal distribution. Hence, by
Remark \ref{rem:4to3geodsfibres}\,, the fibres of $\phi$ are geodesics of $D^M$.
Furthermore, as $(N^3,c_N)$ can be, locally, identified with any leaf
of the horizontal distribution of $\phi$ endowed with the conformal structure induced by $c_M$\,, we have that $(N^3,c_N)$ is flat. Hence, locally, $D^N$ is the Levi-Civita
connection of some constant curvature representative of $c_N$\,. Thus,
if $(M^4,c_M)$ is nonorientable then, up to conformal coverings, $\phi$ is a
harmonic morphisms with geodesic fibres and integrable horizontal distribution
between Riemannian manifolds of constant curvature.
\end{proof}

\begin{exm} \label{exm:GibHaw,Bel,Hit}
1) The harmonic morphisms given by the Gibbons-Hawking and the Beltrami fields
constructions (see \cite{PanWoo-exm}\,) satisfy assertion (i) of Theorem \ref{thm:4to3}\,.\\
\indent
2) The harmonic morphisms of warped product type with one-dimensional fibres from an oriented
four-dimensional Riemannian manifold with nonzero constant sectional curvature satisfy
assertion (ii) of Theorem \ref{thm:4to3}\,. More generally, \emph{let $(M^4,g)$ be the
$\mathcal{H}$-space \cite{LeB-Hspace} of the three-dimensional conformal manifold
$(N^3,c_N)$ and let $\nabla^g$ be the Levi-Civita connection of $g$\,. Suppose that
$(N^3,c_N)$ is endowed with an Einstein-Weyl connection $D^N$ and let $\phi:M^4\to N^3$
be the (local) retraction \cite{Hit-complexmfds} of $N^3\hookrightarrow M^4$ corresponding
to $D^N$\,. Then the map $\phi:(M^4,[g],\nabla^g)\to(N^3,c_N,D^N)$ is a harmonic morphism}
which satisfies assertion (ii) of Theorem \ref{thm:4to3}\,.
\end{exm}

\end{document}